\documentclass[a4paper,12pt]{article} 
\usepackage{amsmath, amsthm, amssymb} 
\usepackage{graphicx,subfigure}
\usepackage[geometry]{ifsym} 
\usepackage{anysize} 
\usepackage{float} 
\usepackage{caption}
\usepackage{amsfonts,amssymb}
\usepackage{epsfig}
\usepackage{verbatim}
\usepackage{color}
\usepackage[normalem]{ulem}

\usepackage{booktabs}
\usepackage{siunitx}
\marginsize{2cm}{2cm}{2cm}{2cm}

\def \d {{\bf d}}

\def \h {{\bf h}}

\def\RR{{\mathop{{\rm I}\kern-.2em{\rm R}}\nolimits}}
\def\PP{{\mathop{{\rm I}\kern-.2em{\rm P}}\nolimits}}

\def\NN{{\mathop{{\rm I}\kern-.2em{\rm N}}\nolimits}}

\def\bfca{\mbox{\boldmath$\alpha$}}
\def\bfcb{\mbox{\boldmath$\beta$}}
\def\bpi{\mbox{\boldmath$\pi$}}

\def\bfmu{\mbox{\boldmath$\lambda$}}

\def \e {{\bf e}}
\newcommand{\be}{\begin{equation}}
\newcommand{\ee}{\end{equation}}
\newcommand{\ba}{\begin{eqnarray}}
\newcommand{\ea}{\end{eqnarray}}
\newcommand{\bi}{\begin{itemize}}
\newcommand{\ei}{\end{itemize}}

\newcommand{\supp}{\mathop{\mathrm{supp}}}

\newcommand{\myspan}{\mathop{\mathrm{span}}}
\newcommand{\trunc}{\mathop{\mathrm{trunc}}}
\newcommand{\myvert}[1]{\!\left.\rule{0ex}{2ex}\right|_{#1}}

\newcommand{\diam}{\mathop{\mathrm{diam}}}

\newtheorem{dfn}{Definition}
\newtheorem{rmk}{Remark}
\newtheorem{lma}{Lemma}

\newtheorem{thm}{Theorem}

\newcommand{\prf}{\noindent{{\bf Proof} :\ }}
\newcommand{\QED}{\vrule height 1.4ex width 1.0ex depth -.1ex\ \medskip}

\def\ACM{Adv.\ Comp.\ Math.}
\def\ACMTG{ACM Trans.\ Graphics}

\def\BIT{BIT}%BIT Numerical Mathematics

\def\CAGD{Comput.\ Aided Geom.\ Design}

\def\CMAME{Comput.\ Methods\ Appl.\ Mech.\ Engrg.}

\def\JAT{J.\ Approx.\ Theory}
\def\JCAM{J.\ Comput.\ Appl.\ Math.}

\def\M3AS{Math.\ Models\ Methods\ Appl.\ Sci.}

\def\NM{Numer. Math.}
\def\VC{Visual Comput.}

\begin{document}

\title{Bivariate hierarchical Hermite spline quasi--interpolation}

\author{{Cesare Bracco}$^{a}$, {Carlotta Giannelli}$^a$,\\ {Francesca Mazzia}$^b$ and {Alessandra Sestini}$^{a}$%\footnote{Corresponding author. E-mail: taewan@snu.ac.kr, Phone: +82-2-880-1437, Fax: +82-2-888-9298.}
\\ \\
\small $^a$Dipartimento di Matematica e Informatica ``U. Dini'', Universit\`a di Firenze\\ \small Viale Morgagni 67/A, 50134 Firenze, Italy \\
\small $^b$Dipartimento di Matematica, Universit\`a di Bari\\ \small Via Orabona 4, 70125 Bari, Italy\\
%\small $^c$Department of Mathematics - Dongguk University \\ \small Pil-dong 3-ga, Jung-gu, Seoul 100-715, Korea\\%\small \texttt{cesare.bracco@yahoo.it}\\
%\and
}
\date{}

%\titlerunning{Spaces of generalized splines over T-meshes}
%\authorrunning{Cesare Bracco \and Fabio Roman}
%\institute{C. Bracco \at
%              Department of Mathematics \lq\lq G. Peano\rq\rq - University of Turin\\ \small V. Carlo Alberto 10, Turin 10123, Italy \\
%              %Tel.: +123-45-678910\\
%              %Fax: +123-45-678910\\
%              \email{cesare.bracco@unito.it}           %  \\
%%             \emph{Present address:} of F. Author  %  if needed
%           \and
%           F. Roman \at
%              Department of Mathematics \lq\lq G. Peano\rq\rq - University of Turin\\ \small V. Carlo Alberto 10, Turin 10123, Italy \\
%              %Tel.: +123-45-678910\\
%              %Fax: +123-45-678910\\
%              \email{fabio.roman@unito.it}              
%}

\maketitle

%******************************* BODY ******************************

%%%%%%%%%%%%%%%%%%%%%%%%%%%%%%%%%%%%%%%%%%%%%
%%%%%%%%%%%%%%%%%%%%%%%%%%%%%%%%%%%%%%%%%%%%%
 
%%%%%%%%%%%%%%%%%%%%%%%%%%%%%%%%%%%%%%%%%%%%%
%%%%%%%%%%%%%%%%%%%%%%%%%%%%%%%%%%%%%%%%%%%%

\begin{abstract}
Spline quasi--interpolation (QI) is a general and powerful approach for the construction of low cost and accurate approximations of a given function. 
In order to provide an efficient \emph{adaptive} approximation scheme in the bivariate setting, we consider quasi--interpolation in hierarchical spline spaces. In particular, we study and experiment the features of the hierarchical extension of the tensor--product formulation of the Hermite BS quasi--interpolation scheme.  The convergence properties of this hierarchical operator, suitably defined in terms of truncated hierarchical B--spline bases, are analyzed. 
A selection of numerical examples is presented to compare the performances of the hierarchical and tensor--product versions of the scheme.
\end{abstract}


\let\thefootnote\relax\footnotetext{Email addresses: cesare.bracco@unifi.it (C. Bracco), carlotta.giannelli@unifi.it (C. Giannelli), francesca.mazzia@uniba.it (F. Mazzia), alessandra.sestini@unifi.it (A. Sestini).}

\section{Introduction} \label{Intro}
Whenever exact interpolation is not suited or not strictly necessary, spline quasi--interpolation (QI) is a valuable and widely appreciated approximation methodology, since QI schemes share \emph{locality} as common denominator  (see e.g., the fundamental papers \cite{debfix,ls75}). 
This is for instance the case when real--time processing  of large data streams is required and the input information on a certain target function $f$ can be dynamically updated (and, consequently, also the spline approximation is dynamically updated).
In  a very general  formulation, a spline quasi--interpolant  $Q(f)$ to a given function $f$ can be written as follows,
\begin{equation} \label{spQIdef}
  Q (f)(\cdot) \,=\,  \sum_{\gamma \in \Gamma} \lambda_{\gamma} (f)\, B_{\gamma}(\cdot),
\end{equation}
where $B_{\gamma}, \,\gamma \in \Gamma\,,$ denote the B--spline basis (or a suitable extension to the multivariate setting) and $\lambda_{\gamma} \,, \gamma \in \Gamma\,,$ are corresponding linear functionals. Even if different approaches can be considered for defining these functionals, they are usually required to be local, namely, 
\begin{equation}  \label{local}
\Lambda_\gamma \,:=\, \supp(\lambda_\gamma) \, \subset \,  \Omega\,, \quad \forall\,\gamma\in\Gamma\,,
\end{equation}  
where $\Omega$ denotes the specific domain of interest.
 The functionals $\lambda_\gamma(f)$ are often defined either as a linear combination of evaluations of $f$ at suitable points or as an appropriate integral of $f$, so that the QI operator exactly reproduces all the polynomials of a certain degree or it is a projector on the space spanned by the B-spline basis (see, e.g., \cite{drs13,llm2001,ls75,sab05a}). QIs where the definition of $\lambda_\gamma(f)$ also involves  derivatives of $f$ are usually called \emph{Hermite QIs} (see, e.g., \cite{MS09}).

QI spline schemes capable to deal with data  on nonuniform meshes are of great interest because they can be profitably used either for adaptive approximation or data reduction purposes. In the univariate setting, the BS Hermite Quasi--Interpolation operator  introduced in \cite{MS09} defines a robust and accurate scheme with this capability.  Its  recent tensor--product extension to the bivariate setting \cite{tesiIurino,ACM15} preserves the robustness and accuracy of the scheme, but necessarily inherits the grid limits of any tensor--product approach. In order to overcome this problem, suitable extensions of tensor--product spline spaces must be considered.
%here we are interested in defining a scheme capable to produce a quasi--interpolating approximation belonging to a hierarchical spline space whose knowledge has been recently deepen in the literature [ ]..... 

Examples of adaptive spline constructions that provide local refinement possibility by relaxing the rigidity of tensor--product configurations have been widely investigated, see e.g.,  \cite{dlp13,gj13,kraft97,ldc10,scfnzl04,sw12}.
%include T--splines \cite{scfnzl04,szbn03}, polynomial splines over T--meshes \cite{ldc10,sw12}, LR--splines \cite{dlp13}, and hierarchical splines \cite{gj13,kraft97}.
The definition of these classes of adaptive splines pose non-trivial issues related to solid basis constructions, the characterization and understanding of the corresponding spline spaces, as well as the design of reliable and efficient refinement algorithms. The introduction of the truncated basis for hierarchical splines has recently provided a powerful construction that can be successfully exploited for spline--based adaptivity and related application algorithms \cite{gjs12,gjs14}.

Among other properties, truncated hierarchical B--splines (THB--splines) preserve  the coefficients of a function represented in terms of the B--spline basis related to a certain hierarchical level \cite{gjs14}. This distinguishing feature can be directly exploited for the definition of hierarchical quasi--interpolants \cite{SM15}. By following this general construction, we study the extension to  hierarchical spline spaces of the above mentioned BS QI scheme. 
The adaptive framework is considered for a hierarchy of nested spline spaces obtained with successive dyadic refinements
and leads to computational benefits with respect to the tensor--product formulation that prevents local mesh refinement.

The paper is organized as follows. Section~\ref{sec_BS} introduces the univariate and tensor--product formulation of the BS quasi--interpolation scheme. A brief overview of hierarchical spline spaces and the THB-spline basis is then introduced in Section \ref{sec_HS}. Section \ref{sec_BSH} is devoted to the hierarchical spline extension of the BS QI scheme and its convergence properties. A selection of numerical experiments, where we employ an automatic refinement strategy to obtain suitable hierarchical meshes, is discussed in Section~\ref{sec_T}. The performances of the BS hierarchical scheme are compared both with its tensor--product counterpart and with a different hierarchical scheme. Finally,  Section \ref{sec_C} concludes the paper.

%%%%%%%%%%%%%%%%%%%%%%%%%%%%

\section{The BS Hermite QI scheme} \label{sec_BS}
 In this section we briefly summarize the Hermite spline QI scheme firstly introduced in \cite{MS09} for the univariate case and recently extended to the bivariate setting in \cite{tesiIurino,ACM15}. A simplified formulation of the scheme is here adopted in order to facilitate the successive introduction of its bivariate hierarchical extension.

\subsection{The univariate scheme}  
In the univariate setting the BS QI scheme approximates any $f:  \Omega = [a\,,\,b] \rightarrow \RR,$ with $f \in C^1(\Omega),$ by defining a suitable spline in  the space  $ S_{d,\pi}, $ of the $d$--degree, $d \ge 2\,,$ $C^{d-1}(\Omega)$  piecewise polynomial functions with breakpoints at the abscissae of a  given  partition  $\pi := \{ x_i, i=0,\ldots,N \}$ of the interval $\Omega.$ The information used are the values of $f$ and $f'$ at the abscissae in $\pi\,.$ Being sufficient for the developed hierarchical generalization, we just relate here to the case of a uniform partition, i.e. to the case $x_i - x_{i-1} = h := (b-a)/N\,, i=1,\ldots,N\,.$ In order to reduce technical difficulties in the next hierarchical formulation, the scheme is further simplified by avoiding the special definition of the first and last $d-1$ functionals which was required in \cite{MS09} (at the negligible price of requiring the knowledge of $f$ and $f'$ at additional $2(d-1)$ abscissae). For this aim, in the following we assume that $f \in C^1( \Omega_e)\,,$ where $ \Omega_e $ is the enlarged interval $ \Omega_e = [a-dh\,,\,b+dh]\,,$ and we require the  knowledge of $f$ and $f'$  at the inner abscissae of the uniform $h$--step  partition  $\pi^e$ of $ \Omega_e, \pi^e = \{x_i\,,  i=-d\ldots,N+d \} \,.$   
Denoting with $B_d$ the $d$--degree B-spline with integer active knots $0,\ldots,d+1\,,\,  Q (f) $ can then be written as follows,
\begin{equation} \label{spvdef}
 Q (f)(\cdot)  \,=\, \bfmu_{d,h}  (f)^T \, {\bf B}_{d,h}(\cdot ) \,=\, \sum_{j=-d}^{N-1} \lambda_j(f)\, B_d(\cdot/h - j),
\end{equation}
where $ \bfmu_{d,h}(f) :=\left(\lambda_{-d} (f),\cdots,\lambda_{N-1} (f) \right)^T\,,$ and where the basis vector ${\bf B}_{d,h}(\cdot) := (B_d(\cdot/h + d),\cdots,B_d(\cdot/h-N+1))^T\,$ selects the elements of the B--spline basis of $S_{d,\pi^e}$ active in $I\,.$  
Note that  the support $\Lambda_j$ of the functional $\lambda_j$ is the following, 
\begin{equation} \label{sup1D}
\Lambda_j  = [x_{j+1}\,,\,x_{j+d}]\,,
\end{equation}  
and that, in the simplified uniform formulation here considered, each functional $ \lambda_j(f)\,, j=-d,\ldots,N-1,$ is defined as the same linear combination of the values $f(x_i)$ and $f'(x_i), i=j+1,\ldots,j+d\,.$ By following the notation used in \cite{MS09}, we can shortly define $\bfmu_{d,h} (f)$ as follows,
\begin{equation} \label{funcBSQ}
\bfmu_{d,h} (f) \,:=\, \hat A ^{(d) } \left( \begin{array}{c}f(x_{-d+1}) \cr \vdots \cr f(x_{N+d-1}) \end{array}  \right)\,-\, 
h\, \hat B ^{(d)} \left( \begin{array}{c} f'(x_{-d+1}) \cr \vdots \cr f'(x_{N+d-1}) \end{array} \right)\,,
\end{equation}
where $ \hat A ^{(d)}\,, \hat B ^{(d)} \in \RR^{(N+d) \times (N+2d-1)}$ are two Toeplitz banded matrices of bandwidth $d$ such that
$$ \e_j^T\, \hat A ^{(d)} =  ({\bf 0}^T_{j-1}\,, \hat{\bfca}^{ (d)T } \,, {\bf 0}^T_{N+d-j})\,, \qquad \e_j^T\, \hat B ^{(d)} =  ({\bf 0}^T_{j-1}\,, \hat{\bfcb}^{ (d)T } \,, {\bf 0}^T_{N+d-j})\,,$$ 
with $(\hat{\bfca}^{(d)T} \,,  \hat{\bfcb}^{(d)T})^T \in \RR^{2d}$ obtained by solving a suitable local nonsingular $2d \times 2d$ linear system.
\begin{table}
\label{tabcoef}
\centering
\begin{tabular} {c|c c}
 $d$  & $\hat{\bfca}^{(d)T}$  &$\hat{\bfcb}^{(d)T}$ \\
\hline
 $2$ & $\frac{1}{2}\,(1\,,\,1)$ &   $\frac{1}{4}\,(-1\,,\,1)$ \\
 $3$ & $\frac{1}{2}\,(-1\,,4\,,-1)$ &   $\frac{1}{6}\,(1\,,0\,,-1)$ \\
 $4$ & $\frac{1}{12}\,(5\,,1\,,1\,,5)$ &   $\frac{1}{48}\,(-5\,,-41\,,41\,,5)$ \\
\end{tabular}
\caption{Vector coefficients defining the uniform $d$--degree BS Hermite spline quasi--interpolant.}
\end{table}
We relate to \cite{MS09} for the details about such system and we just report in Table 1 the values of the vector coefficients $\hat{\bfca}^{(d)}$ and $\hat{\bfcb}^{(d)}$ for $2 \le d \le 4\,.$ Note that in \cite{MS09} it has been proved that the operator $ Q_{d}$ is a projector  which means that $ Q_{d\,,\,\pi}(s) = s\,, \, \forall s \in S_{d,\pi}$ and so, in particular  $ Q_{d} (p) = p\,, \, \forall p \in \PP^d\,.$  
Concerning the convergence, the analysis developed in \cite{MS09} under the aforementioned hypotheses can be simplified to $|\lambda_j(f)| \,\le \,  \| \hat{\bfca}^{(d)} \|_1 \, \|f\|_{\Lambda_j} + h   \| \hat{\bfcb}^{(d)} \|_1\, \|f'\|_{\Lambda_j}\,,$ since
$$ \|\hat A^{(d)}\|_{\infty} = \| \hat{\bfca}^{(d)}  \|_1 \,, \qquad  \|\hat B^{(d)}\|_{\infty} = \| \hat{\bfcb}^{(d)}  \|_1\,,$$
where, 
\[
\|g\|_{c} := \max_{u \in c} |g(u)| \,, \quad  \forall  g \in C^0(\Omega)\,. 
\]
 Considering the locality, the non negativity and the partition of unity property fulfilled in $\Omega$ by the B--splines selected in ${\bf B}_d\,,$ this result implies also that for any cell $c =[x_i\,,\,x_{i+1}], i=0,\ldots,N-1,$ and for any function $f \in C^1(\Omega)$ it is
$$\|Q(f)\|_{c} \le   \| \hat{\bfca}^{(d)} \|_1 \, \|f\|_{C} + h   \| \hat{\bfcb}^{(d)} \|_1\, \|f'\|_{C}\,,$$
where $C = \displaystyle{\bigcup_{j=i-d}^i} \Lambda_i = [x_{i-d+1}\,,\,x_{i+d}]\,,$ whose size  $|C|$ is equal to $(2d-1)h\,.$ 

By writing $\|f -  Q(f) \|_{c} \le  \|f -  p \|_{c} + \| Q(p-f) \|_{c} \quad \forall p \in \PP^d\,,$   we can conclude that
$$\|f -  Q(f) \|_{c} \le  (1 + \| \hat{\bfca}^{(d)} \|_1)\,  \|f -  p \|_{C} \,+\, h   \| \hat{\bfcb}^{(d)} \|_1\, \|f' - p'\|_{C}\,. $$
Now, by selecting $ p \in \PP^d$ as the Taylor expansion of $f$ at some point in $c,$ considering that $|C| = (2d-1)h,$ we have that   $\| f -  p \|_{C} =  \mathcal{O}(h^{k+1})$ and $\|f' - p'\|_{C} =  \mathcal{O}(h^d)\,.$ Since $c$ can be any cell in $\Omega,$ we can conclude that the scheme has maximal approximation order, that is the error $\|f -  Q(f) \|_{\Omega}$ is $ \mathcal{O}(h^{d+1})\,,$  when $f \in C^{d+1}( \Omega_e).$

%----------------------------------

\subsection{The bivariate scheme}
When $f$  is a sufficiently smooth bivariate function defined on a rectangle, $f=f(x,y)$ with $f: \Omega := I_1 \times I_2 \rightarrow \RR\,$ with $I_j = [a_j\,,\,b_j], j=1,2,$   we can apply the univariate BS quasi--interpolation operator in  any of the two directions, thus obtaining two different approximations of $f.$ If the scheme is applied to $f$ twice, once with respect to $x$ and the other to $y\,,$ the BS tensor--product spline quasi--interpolant is obtained. Clearly, this bivariate extension of the scheme needs of two partitions, $\pi_1 = \{x_0,\ldots,x_{N_1} \}$ of $I_1$ and   $\pi_2 = \{y_0,\ldots,y_{N_2} \}$   of   $I_2\,,$ specifying the $x$ and $y$ spline breakpoints,  besides of two corresponding integers, $d_1$ and $d_2\,,$ prescribing the spline degree with respect to $x$ and to $y.$ For the more general formulation of the bivariate BS QI scheme we refer to \cite{tesiIurino,ACM15}. Here we briefly summarize it, specifically relying on  the uniform univariate formulation of the method described in the previous section  and consequently we have $N_1 \,:=\, (b_1-a_1)/h_x\,, \,\, N_2 \,:=\, (b_2-a_2)/h_y\,.$
Then we actually need two  enlarged partitions  $\pi_1^e =  \{x_{-d_1},\ldots,x_{N_1+d_1} \} \supset \pi_1$ and   $\pi_2^e =  \{y_{-d_2},\ldots,y_{N_2+d_2} \} \supset \pi_2 \,$ of the two enlarged intervals $I_1^e \supset I_1$ and $I_2^e \supset I_2\,,$ and we require $f \in C^{(1,1)}( \Omega_e)\,,$ where  now $ \Omega_e \,:=\, I_1^e \times I_2^e\,.$  By using the notation introduced in the previous section,  the two intermediate approximations $Q_1  (f)$ and  $Q_2  (f),$ can be written as follows,  
$$\begin{array}{ll} Q_1  (f) (x,y) & =  [\bfmu ^1_{d_1,h_x}  (f)] (y) ^T \, {\bf B}_{d_1,h_x} (x)\,, \cr
Q_2  (f) (x,y) & =  [\bfmu^2_{d_2,h_y}  (f)] (x) ^T \, {\bf B}_{d_2,h_y} (y) =  {\bf B}_{d_2,h_y} (y)^T \,  [\bfmu^2_{d_2,h_y}  (f)] (x)\,,
 \cr
\end{array}
$$
where we have emphasized that  $Q_1 (f)$ ($Q_2 (f)$) is a spline with respect to $x$ ($y$) and a general function with respect  to the other variable. 
Thus the final tensor--product  approximation $   (Q_2  Q_1) (f) $ is obtained by applying the operator $Q_2 $ to the intermediate approximation $Q_1(f)$ (or alternatively doing the opposite) and clearly it belongs to the tensor--product spline space $ S_{\d,\h}\,,$ where
$$S_{\d,\h} \,:=\, S_{d_1, \pi_1} \otimes S_{d_2,\pi_2}\,, \mbox{ with } \d \,:=\, (d_1\,,\,d_2)\,, \,\,\h \,:=\, (h_x\,,\,h_y)\,.$$
Following this procedure, with some computations it can be verified that, denoting just with $Q$ the operator $Q_2 Q_1\,,$  $ Q (f)$ can be compactly  written as follows,
\begin{equation} \label{TP}   Q (f)(x,y)   \,=\, {\bf B}_{d_2,h_y} (y)^T \,M_{\d,\h}(f)\,\,  {\bf B}_{d_1,h_x} (x) \,,
\end{equation}
where   
$$  \begin{array}{ll}
M_{\d,\h}(f)  \,:=\,   &\hat A ^{(d_2)} \,  F \, \hat A ^{(d_1)T}\   
\,-\, h_x\,  \hat A ^{(d_2)}  \, F_x \, \hat B ^{(d_1)T}  \,-\, \cr
 &h_y\, \hat B ^{(d_2)}  \, F_y \, \hat A^{(d_1)T} 
\,+\,   h_x\,h_y \, \hat B ^{(d_2)}  \, F_{xy} \, \hat B^{(d_1)T} \,.
\end{array}  $$
The matrices $F, F_x, F_y$ and $F_{xy}$ involved in the above formula are full matrices belonging to $\RR^{(N_2+2d_2-1) \times (N_1+d_1-1)}$ such that,
$$
F \,:=\,   \left( \begin{array}{ccc} f(x_{-d_1+1},y_{-d_2+1}) & \cdots & f(x_{N_1+d_1-1},y_{-d_2+1}) \cr
\vdots & \ddots & \vdots \cr  f(x_{-d_1+1},y_{N_2+d_2-1}) & \cdots & f(x_{N_1+d_1-1},y_{N_2+d_2-1}) \end{array} \right) \,,
$$
and $F_x, F_y$ and $F_{xy}$ are analogously defined just replacing $f$ respectively with $f_x, \,\,f_y$ and with $f_{xy}\,.$
Thus the considered tensor--product formulation of the  BS QI scheme needs the knowledge of $f, f_x, f_y, f_{xy}$ in the inner points of the extended lattice $\bpi^e := \pi_1^e \otimes \pi_2^e\,.$   In order to facilitate the introduction of the hierarchical generalization of our scheme, we are interested in rewriting  the expression of $Q(f)$ defined in  (\ref{TP}) in the following  form,  
\begin{equation} \label{TPalt}   Q (f) (x,y)   \,=\,  \sum_{J \in \Gamma_{\d,\h}} \lambda_J(f) \, B_J^{(\d,\h)}(x,y)\,,
\end{equation}
where the basis function $B_J^{(\d,\h)}$ associated to the multi--index $J =(j,i)$ is the following  product of univariate B-spline functions,
\begin{equation} \label{TPbasis}
B_J^{(\d,\h)}(x,y) \,:=\, B_{d_1}(x/h_x-i)\,\,  B_{d_2}(y/h_y-j)\,,
\end{equation}
and where  
\begin{equation} \label{Gammadef} \Gamma _{\d,\h}\,:=\, \{ (j,i), j = -d_2+1,\ldots,N_2+d_2-1\,,\,\,  i = -d_1+1,\ldots,N_1+d_1-1 \}\,.\end{equation}
 
From  (\ref{TP}), recalling the structure of the matrices $\hat A_{d_i}, \hat B_{d_i}, i=1,2\,,$ we obtain the following definition for the functional $\lambda_J(f)$ associated with the multi--index $J =(j\,,\,i)\,,$

\begin{equation} \label{funcdef}
 \begin{array}{ll} \lambda_J(f) \,=\,
&\hat \bfca^{(d_2)T}\,  F^J \, \hat \bfca^{(d_1)} -   h_x \, \hat \bfca^{(d_2)T}\, F^J_x\, \hat \bfcb^{(d_1)} - \cr
& h_y \, \hat \bfcb^{(d_2)T}\,  \, F^J_y\, \hat \bfca^{(d_1)} +  h_x h_y \, \hat \bfcb^{(d_2)T}\,  F^J_{xy}\, \hat \bfcb^{(d_1)} \,, \cr
\end{array} 
\end{equation}
where the local matrix $F^J \in \RR^{d_2 \times d_1}$ is defined as follows,
$$
F^J \,:=\,   \left( \begin{array}{ccc} f(x_{i+1},y_{j+1}) & \cdots & f(x_{i+d_1},y_{j+1}) \cr
\vdots & \ddots & \vdots \cr  f(x_{i+1},y_{j+d_2}) & \cdots & f(x_{i+d_1},y_{j+d_2}) \end{array} \right) \,,
$$
and where $F^J_x, F^J_y$ and $F^J_{xy}$ are analogously defined just replacing $f$ respectively with $f_x, \,\,f_y$ and with $f_{xy}\,.$   This implies that, using the notation introduced in (\ref{local}), we have that  the area of $\Lambda_J$, support of $\lambda_J$, is equal to   $(d_1-1) (d_2-1) \,h_xh_y\,,$ as it is  
\begin{equation} \label{suppfunc}
\Lambda_J \,=\, [x_{i+1}\,,\,x_{i+d_1}] \times [y_{j+1}\,,\,y_{j+d_2}]\,.
\end{equation}
The following lemma gives an upper bound for the absolute value of each local linear functional $\lambda_J(f)$ and it will be later useful to develop the study of the approximation power of the hierarchical extension of $Q(f),$
\begin{lma} \label{lemmu}
There exist four positive constants $\kappa_{r,s}, r\, ,s = 0,1$ depending on $\d$ but not depending on $\h\,,$ such that for any function $f \in C^{(1,1)}(\Omega_e),$ and for each multi--index $J \in \Gamma_{\d,\h}$ it is
\begin{equation} \label{muineq}
|\lambda_J(f)| \le \kappa_{0,0} \|f\|_{\Lambda_J} +\, h_x\, \kappa_{1,0} \|f_x\|_{\Lambda_J} +\, h_y\, \kappa_{0,1} \|f_y\|_{\Lambda_J} +\,    h_x h_y \kappa_{1,1} \|f_{xy}\|_{\Lambda_J} \,. 
\end{equation}
\end{lma}
\prf
Clearly we can upper bound $|\lambda_J(f)|$ with the sum of the absolute values of the four addends defining it in (\ref{funcdef}). Now, relating for example to  the first of these addends, we have that
$$ | \hat \bfca^{(d_2)T}\,  F^J \, \hat \bfca^{(d_1)} | \,\le\, \kappa_{0,0} \| f \|_{\Lambda_J}$$
where $\kappa_{0,0} := \| \hat \bfca^{(d_2)} \|_1\,  \| \hat \bfca^{(d_1)} \|_1 \,.  $
Using similar arguments for the other three addends the proof can be easily completed.
\QED
 
Note that, being  the univariate BS QI operator a projector, also the corresponding tensor--product operator is still a projector which means that $Q(s) = Q_2Q_1(s) = s $ when $s \in S_{\d,\h} \,.$ Concerning the approximation order, when $f \in C^{(d_1+1,d_2+1)}( \Omega_e)\,,$  it can be proved that the error $\|  Q (f) - f \|_{\Omega}$ is $ \mathcal{O}(h_x^{d_1+1} ) \,+\, \mathcal{O}( h_y^{d_2+1} )\,.$   Actually, denoting with $I$ the identity operator, this can be easily verified considering that, for any tensor--product operator $P = P_1P_2\,,$ the corresponding  error operator $ E:= I - P\,,$ can be split as follows, 
\begin{equation} \label{split} E \,=\, E_1 + E_2 -  E_1E_2\,,
\end{equation}
where $E_i := I -P_i ,\,\, i=1,2\,.$  

 \begin{rmk}
The above adopted assumption to know $f$ on the enlarged domain $ \Omega_e$ is not essential to approximate $f$ in $\Omega$ by using the BS QI scheme either for the monovariate or the bivariate case \cite{MS09,tesiIurino,ACM15}. It is here assumed just to simplify the introduction of the later developed bivariate herarchical extension of the scheme.
\end{rmk}

%%%%%%%%%%%%%%%%%%%%%%%%%%%%%%%%%%%%
\section{The bivariate hierarchical space}  \label{sec_HS}
Hierarchical bivariate B-spline spaces are constructed by considering a hierarchy of $M$ tensor--product bivariate spline spaces $V^0\subset V^1\subset \ldots...\subset V^{M-1}$ defined on $\Omega$, together with a hierarchy of closed domains $\Omega^0 \supseteq \Omega^1 \supseteq \dots\supseteq\ \Omega^{M-1}$, that are subsets of $\Omega,$ with $\Omega^0 =\Omega.$ The \emph{depth} of the subdomain hierarchy is represented by the integer $M$, and we assume $\Omega^M=\emptyset$. Note that each $\Omega^{\ell}$ is not necessarily a rectangle, and it may also have several distinct connected components. Even if  hierarchical spline constructions can cover different non-uniform knot configurations and degrees according to the considered nested sequence of spline spaces, we restrict our analysis to dyadic (uniform) refinement and a fixed degree $\d$ at any level, a typical choice in hierarchical approximation algorithms.

\begin{figure} 
\centering\includegraphics[trim=0mm 0mm 0mm 0mm,clip,scale=0.40]{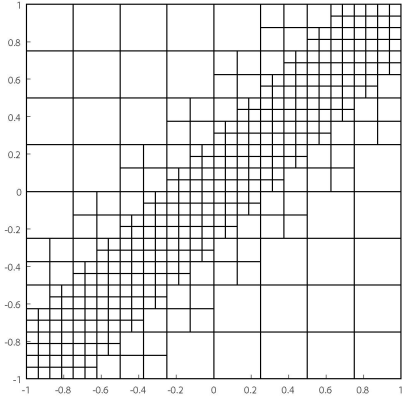}   %change figure
\caption{A hierarchical mesh obtained with $M=3$ levels.}
 \label{Hiermesh}
 \end{figure}

The spline spaces $V^\ell$, for $\ell=0,\ldots,M-1$, are defined through successive dyadic refinements of an initial uniform tensor--product grid characterized by cells of size ${h_x}\times{h_y}$, namely
$$
V^{\ell} := S_{\d,\h_{\ell}}\,, \quad 
\h_{\ell} :=\left(h_{x,\ell}\,,\,h_{y,\ell}\right) := \left(\frac{h_x}{2^{\ell}}\,,\frac{h_y}{2^{\ell}}\right)\,.
$$
Thus, for each level $\ell$,   $V^\ell$ is spanned by  the tensor--product B-spline basis,
$${\cal B}_\d^{\ell} \,:=\, \{ B_J^{ \ell }, J \in  \Gamma_{\d,\h_{\ell}} \}\,, $$
where $ B_J^{\ell}$ simply denotes the function $B_J^{(\d,\h_{\ell})}$ defined in (\ref{TPbasis}), and $ \Gamma_{\d,\h_{\ell}} $ is defined in (\ref{Gammadef}).
Let $G^\ell$ be the uniform tensor--product grid with respect to level $\ell$ defined as the collection of cells $c$ of area $h_x^\ell\, h_y^\ell$ covering $\Omega$. As we adopt the common assumption that each $\Omega^{\ell}$, $\ell\ge 1$, can be obtained as the union of a finite number of cells  $c$ of the grid $G^{\ell-1}\,,$ then we can define the \emph{hierarchical mesh} ${\cal G}_{\cal H}$ as
\[
{\cal G}_{\cal H} := %\bigcup_{\ell=0,\ldots,M-1}  {\cal G}^\ell
\left\{ c\in {\cal G}^\ell, \forall\,\ell=0,\ldots,M-1\right\}
\]
where
\[
 {\cal G}^\ell := \left\{c\in G^\ell : c\subseteq \Omega^\ell \wedge 
\nexists\, c^*\in G^{\ell^*},\, \ell^*>\ell: c^*\subseteq\Omega^{\ell^*} \,\wedge\,
 c^*\subset c\right\},  
\]
for $\ell=0,\dots,M-1$, are the sets that collect the active cells at different levels. Figure \ref{Hiermesh} shows an example of a sequence of nested domains for $M=3$ and the corresponding hierarchical mesh ${\cal G}_{\cal H}$.

In order to define the spline hierarchy,  the set of active basis functions $B_J^{\ell}$, for each level $\ell$, is selected from ${\cal B}_\d^{\ell}$ according to the following definition, see also \cite{kraft97,vgjs11}.
\begin{dfn}\label{dfn:hb}
The hierarchical B-spline basis ${\cal H}_\d({\cal G}_{\cal H})$ of degree $\d$ 
with respect to the mesh ${\cal G}_{\cal H}$ is defined as
\begin{equation*}
{\cal H}_\d({\cal G}_{\cal H})  :=  \left\{
B_J^{\ell} \in {\cal B}_\d^\ell \, : \, J \in A_\d^{\ell}, \,\, \ell=0,...,M-1  \right\} \,\,,
\end{equation*} 
where  
\begin{equation} \label{acitve}
A_\d^{\ell} \,:=\, \{ J \in \Gamma_{\d,\h_{\ell}}\,:\,  
\supp B_J^{\ell}  \subseteq \Omega^\ell\,\,\wedge \,\, 
\supp B_J^{\ell} \not\subseteq \Omega^{\ell+1} \}\,,
\end{equation}
is the active set of multi--indices $A_\d^{\ell} \subseteq  \Gamma_{\d,\h_{\ell}}$ and $\supp B_J^{\ell}$ denotes the intersection of the support of $B_J^{\ell}$ with $\Omega^0$.
\end{dfn}

The space  $S_{\cal H} := \myspan {\cal H}_\d({\cal G}_{\cal H}) $ is the hierarchical spline space associated to the mesh ${\cal G}_{\cal H}$. Obviously, $S_{\cal H}$ is a subset of  $V^{M-1} \subseteq C^{(d_1-1,d_2-1)}(\Omega)$.  In virtue of the linear independence of hierarchical B--splines, see e.g., \cite{kraft97,vgjs11},
$$
\dim(S_{\cal H} )\,=\,  \sum_{\ell = 0}^{M-1} \# A_\d^\ell\, \, \le \dim (V^{M-1}).
$$
Note that, the symbol $\ll$ can usually replace the last $\le$ inequality symbol, since the adaptive spline framework allows to obtain the desired approximation with a strongly reduced number of degrees of freedom.  Furthermore we would like to remark that, even if in general it is only $S_{\cal H} \subseteq S_{\d}({\cal G_{\cal H}})$, where
\[
S_{\d}({\cal G}_{\cal H})
:= \left\{
f\in C^{(d_1-1,d_2-1)}(\Omega): f\myvert{c}\in\PP^\d, \forall c\in{\cal G}^\ell,
\ell=0,\ldots,M-1
\right\},
\]
under reasonable assumptions on the hierarchical mesh configuration it is possible to obtain that  $S_{\cal H} = S_{\d}({\cal G_{\cal H}})$ \cite{gj13}.

In order to recover the partition of unity property in the hierarchical setting and reduce the support of coarsest basis function, the truncated hierarchical B-spline (THB--spline) basis of $S_{\cal H}$ has been recently introduced \cite{gjs12}. This alternative hierarchical construction relies on the definition of the following truncation operator $\trunc^{\ell+1} : V^{\ell} \rightarrow V^{\ell+1}, \ell=0,\ldots,M-2$.

\begin{dfn} \label{truncoperator}
Let 
$ s(\cdot) = \displaystyle{\sum_{J \in \Gamma_{\d,\h_{\ell+1}} } } \sigma_J^{\ell+1} B_J^{\ell+1}(\cdot)\,,$
be the representation in the B--spline basis of $V^{\ell+1} \supset V^\ell\,,$ of
$s \in V^\ell$. The truncation of $s$ with respect to ${\cal B}^{\ell+1}$ and $\Omega^{\ell+1}$ is defined as 
\begin{equation*}
{\trunc}^{\ell+1} s (\cdot) \,:=\, 
\displaystyle{ \sum_{J \in\Gamma_{\d,\h_{\ell+1}} \,:\, 
\supp B_J^{\ell+1}\, \not\subseteq\, \Omega^{\ell+1}  } } \sigma_J^{\ell+1} B_J^{\ell+1}(\cdot)\,.
\end{equation*}
 \end{dfn}
 
By also introducing the operator $Trunc^{\ell+1} : V^{\ell} \rightarrow S_{\cal H} \subseteq V^{M-1}$, for $\ell=0,\ldots,M-1,$
  $$Trunc^{\ell+1} \,:=\,   trunc^{M-1}( trunc^{M-2}(\cdots( trunc^{\ell+1}(s)) \cdots))\,,$$  
we are ready to define the truncated basis for $S_{\cal H}\,, $ according to the following construction \cite{gjs12}.

\begin{figure} \begin{center}
\subfigure[]{\includegraphics[trim=25mm 5mm 25mm 8mm,clip,scale=0.53]{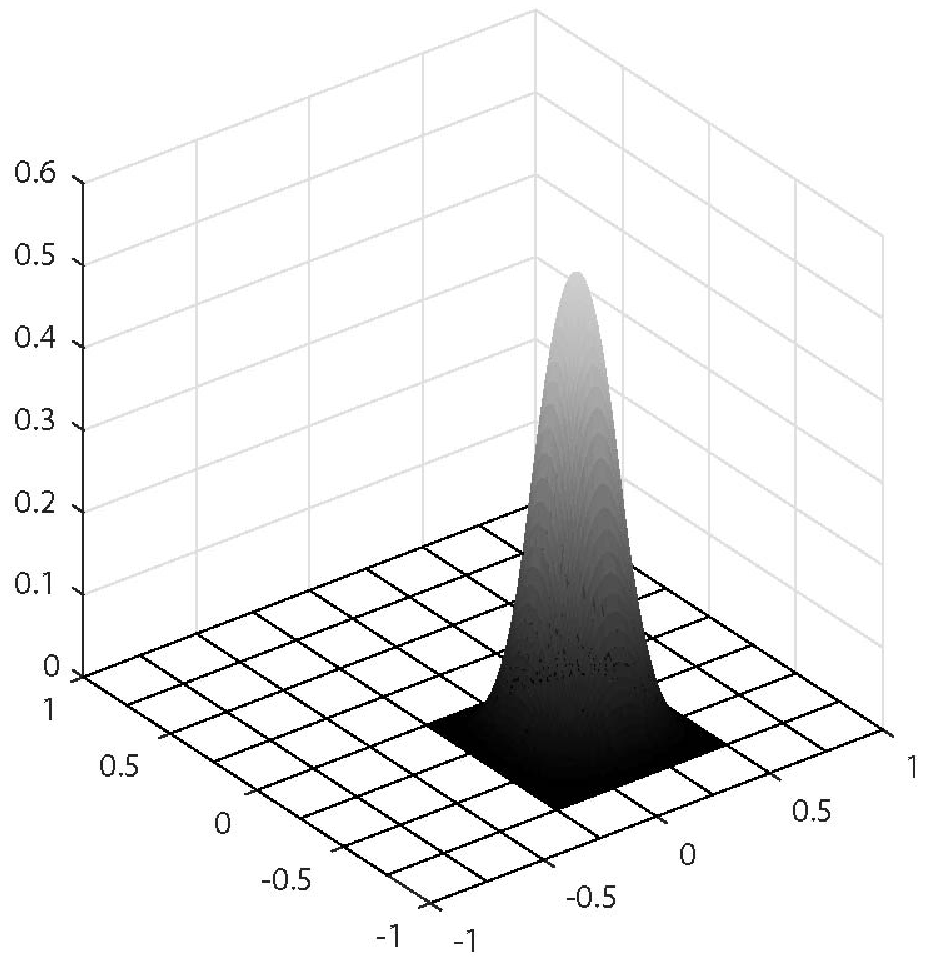}}  %change figure
\subfigure[]{\includegraphics[trim=25mm 5mm 25mm 8mm,clip,scale=0.53]{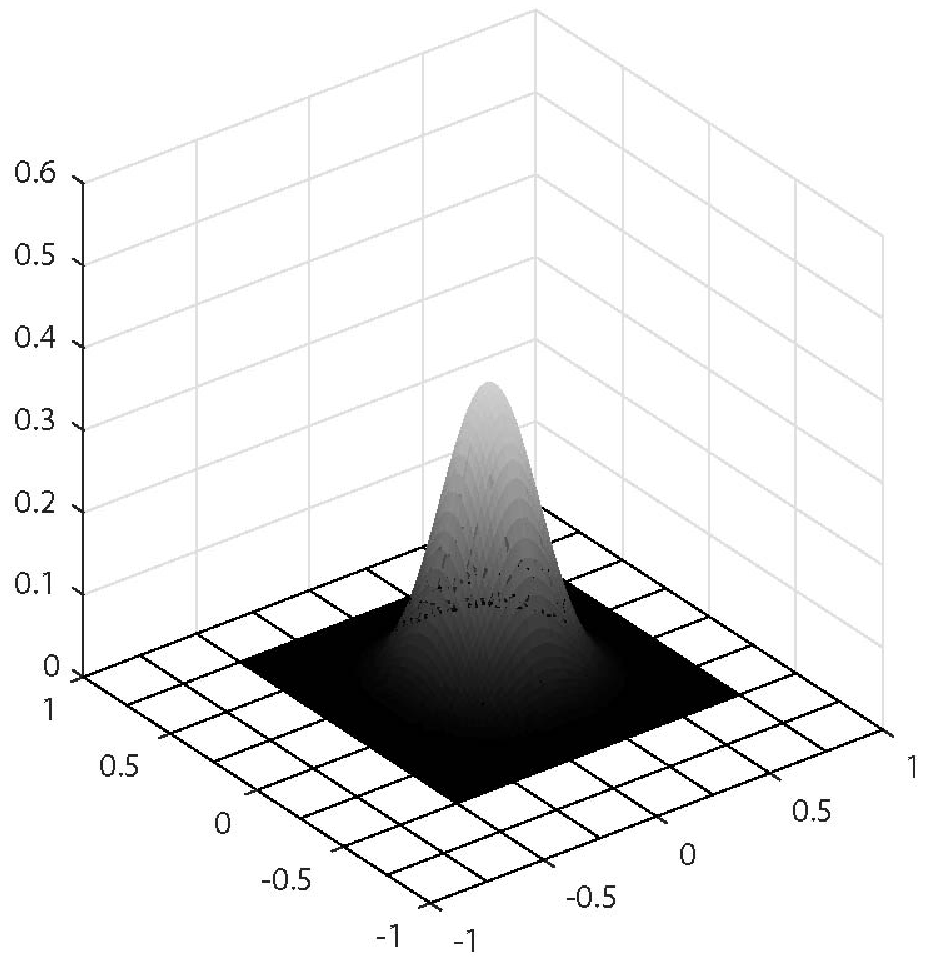}}  %change figure
\\
\subfigure[]{\includegraphics[trim=25mm 5mm 25mm 8mm,clip,scale=0.53]{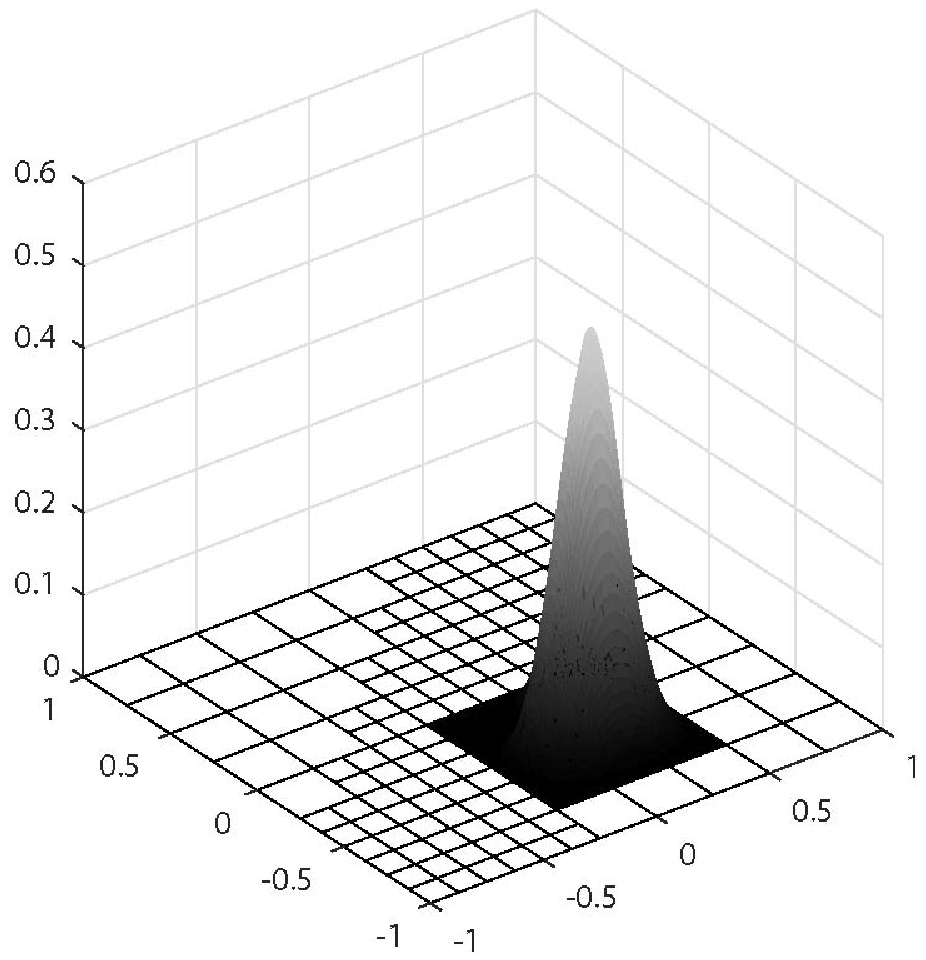}}
\subfigure[]{\includegraphics[trim=25mm 5mm 25mm 8mm,clip,scale=0.53]{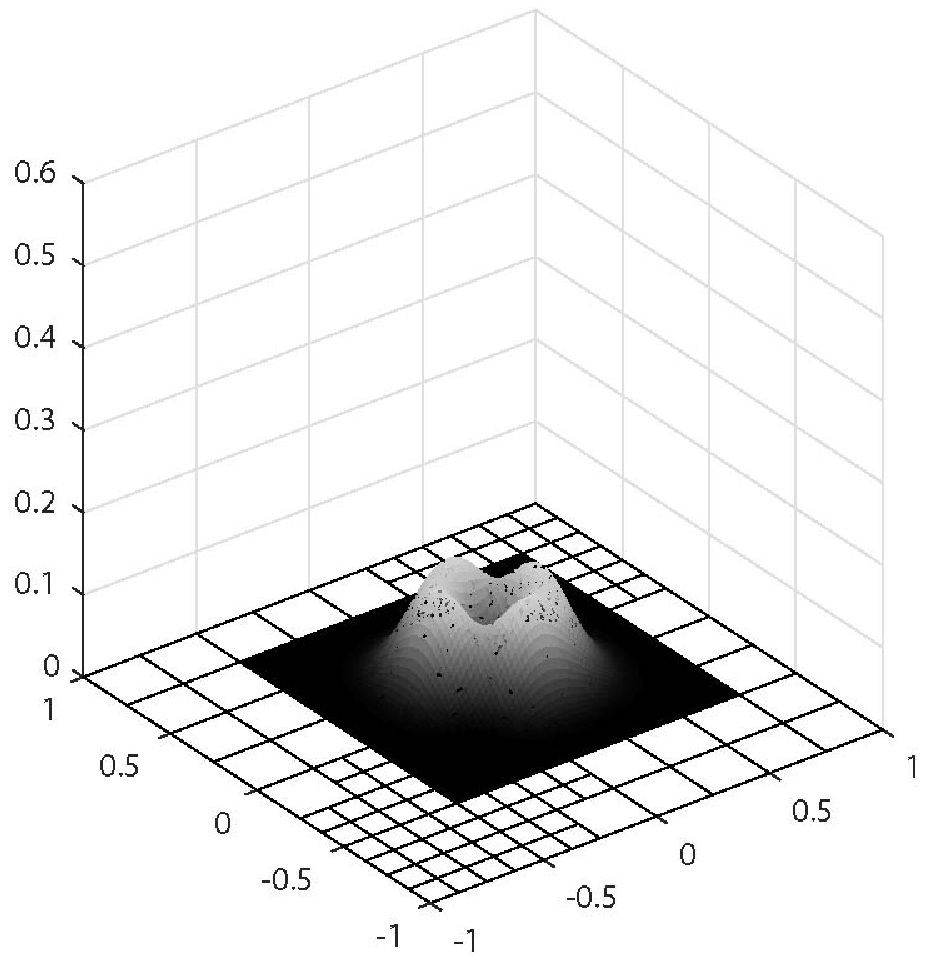}}
\\
\subfigure[]{\includegraphics[trim=25mm 5mm 25mm 8mm,clip,scale=0.53]{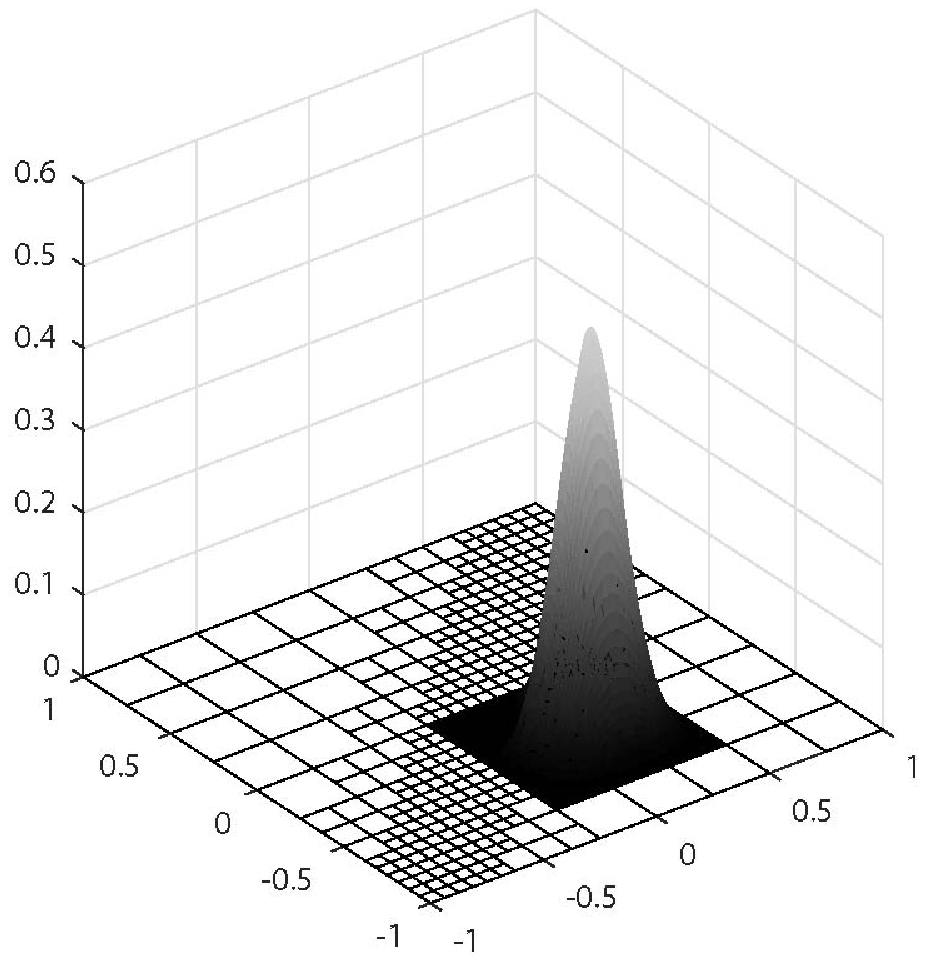}}
\subfigure[]{\includegraphics[trim=25mm 5mm 25mm 8mm,clip,scale=0.53]{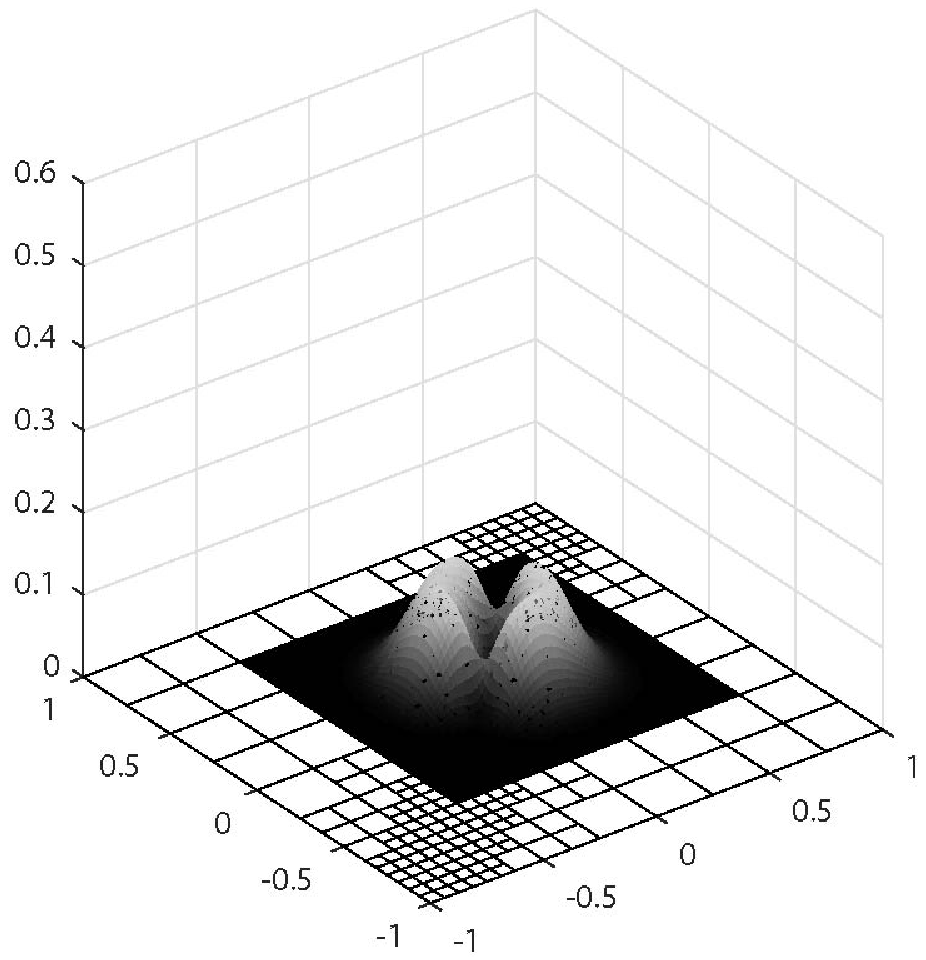}}
\caption{Two examples of elements of the truncated hierarchical B-spline basis. For $(d_1,d_2)=(2,2)$, starting from (a) $B_J^0$, the truncation at level $\ell=1$ gives (c) $trunc^{1}(B_J^0)$, and the further truncation at level $\ell=2$ gives (e) $T_J^{2}:=Trunc^{2}(B_J^0):=trunc^{2}(trunc^{1}(B_J^0))$. Similarly, (b), (d), (f) represent the successive truncations of a B-spline of bi-degree $(d_1,d_2)=(4,4)$.}\label{truncfig}
\end{center}
\end{figure}

\begin{dfn}\label{dfn:thb}
The truncated hierarchical B-spline basis ${\cal T}_{\d} ({\cal G}_{\cal H})$ of degree $\d$  with respect to the mesh ${\cal G}_{\cal H}$ is defined as
\begin{equation}  \label{truncbasis}
{\cal T}_{\d} ({\cal G}_{\cal H}) \,:=\,
\left\{ T_J^{\ell} \,:\, J \in A_\d^{\ell}, \,\, \ell=0,...,M-1  \right\}\,, \,\, \text{with} \quad T_j^{\ell} \,:=\, Trunc^{\ell+1} (B_J^\ell)\,.  
\end{equation}
 \end{dfn}
 THB--splines have very important properties like linear independence, non--negativity, partition of unity, preservation of coefficients and strong stability \cite{gjs12,gjs14}. In particular the  preservation of  coefficients property  is of central interest in the context of quasi--interpolation because it ensures that  the hierarchical extension of any tensor--product quasi--interpolation operator keeps its polynomial reproduction capability \cite{SM15}.
 Figure \ref{truncfig} shows the two successive steps necessary to obtain the truncated basis element $T _J^0$ for a certain $J \in A_\d^0$ with respect to the hierarchy of domains shown in Figure 1.
 Note that, if $\sigma_J^{\ell}$ denotes the coefficient of a spline $s \in S_{\cal H}$  in the THB spline basis, since only few truncated basis functions do not vanish on a cell  $c$ of the hierarchical mesh, we can write,
$$ s(x,y) \,=\, \sum_{\ell = 0}^{M-1}  \sum_{J \in {\cal A}_\d^\ell(c) }\sigma_J^\ell T_J^\ell (x,y)\,, \quad \forall (x,y) \in c\,,$$
 where %\textcolor{red}{dove usiamo questa rappresentazione? ci serve?}  
 \begin{equation} \label{localex}
  {\cal A}_\d^\ell(c) \,:=\, \{J \in A_\d^\ell\,\, |\,\, c \subseteq supp(T_J^\ell) \} \,. 
  \end{equation}

%\vspace{-1.0cm}
 
The truncation of coarsest basis functions with respect to finer basis elements allow to reduce the overlap of supports of basis functions at different hierarchical levels. Nevertheless, in order to have a specific bound for the number of truncated basis functions acting on any cell, suitable hierarchical configurations should be selected.

\begin{dfn}\label{dfn:amesh}
A hierarchical mesh ${\cal G}_{\cal H}$ is admissible of class $m$, $1\le m<M$, if the THB--splines which take non--zero value over any cell $c\in{\cal G}_{\cal H}$ belong to at most $m$ successive levels.
\end{dfn}

Restricted hierarchies that allow the identification of admissible meshes were considered in \cite{gjs14} for the case $m=2$, and in \cite{bg15} for a general $m$.  If ${\cal G}_{\cal H}$ is an admissible mesh of class $m$ the following remarks can be done,
\begin{itemize}
\item if $c \in {\cal G}_k\,,$  ${\cal A}_\d^\ell(c)$ can be non empty only if $k_m \le \ell \le k\,,$ with
\begin{equation} \label{kmdef}
k_m \,:=\, \max\{0\,,\,k-m+1\}\,;
\end{equation} 
\item $m(d_1+1)(d_2+1)$ is an upper bound for the number of THB--splines which are non--zero on each cell; 
\item there exist two positive constants $\gamma_1$ and $\gamma_2$ depending on the class of admissibility $m$ but not on the other quantities such that, for all $c\in{\cal G}_{\cal H}$ and for each multi--index $J \in {\cal A}_\d^\ell(c)\,,$  it is  $$\gamma_1 \diam(c) \le \diam(\supp T_J^\ell) \le \gamma_2\diam(c)\,.$$ 
 
\end{itemize}
The above properties of hierarchical meshes of class $m$ will be used in the next section  to get theoretical proofs of convergence for the considered hierarchical quasi--interpolation scheme.

%%%%%%%%%%%%%%%%%%%%%%%%%%%

\section{The hierarchical Hermite BS QI scheme}  \label{sec_BSH}
Besides the mentioned advantages,  the THB--spline basis of $S_{\cal H}$ has another fundamental advantage with respect to the HB--spline basis  because, thanks to its preservation of coefficients property, it  allows an  easy  extension to the hierarchical setting of any quasi--interpolation operator \cite{SM15}. For example, specifically relating  the tensor--product BS QI operator $Q$ defined in (\ref{TPalt}), we can just define its extension $Q_{\cal H} : C^{1,1}(\Omega_e) \rightarrow S_{\cal H}$ as follows,
\begin{equation} \label{QH}   Q_{\cal H} (f) (x,y)  \,:=\,  \sum_{\ell = 0}^{M-1} \sum_{J \in A_{\d,\h_\ell}} \lambda_J^\ell(f) \, T_J^{\ell} (x,y)\,,
\end{equation}
with $T_J^\ell, \,J \in A_{\d,\h_\ell}\,, \ell=0,\ldots,M-1,$ denoting the truncated basis introduced in (\ref{truncbasis}) and with the functional $\lambda_J^\ell(f)$ defined as in (\ref{funcdef}), just replacing $h_x$ and $h_y$ respectively with $h_{x,\ell}$ and $h_{y,\ell}\,.$
Actually  in \cite{SM15} it has been proved  that this definition is reasonable because it ensures the following implication,
 \begin{equation} \label{rippol}
 Q(p) \,=\, p\,,\,\, \forall p \in \PP^\d \quad \Rightarrow \quad Q_{\cal H}(p) = p\,,\,\, \forall p \in \PP^\d\,.
 \end{equation} 
 We start with a preliminary lemma.

\begin{lma} \label{lem2}
If the hierarchical mesh ${\cal G}_{\cal H}$ is admissible of class $m\,,$ then  for any cell $c \in {\cal G}^k, k=0,\ldots,M-1,$ and for any function $f \in C^{(1,1)}(\Omega_e),$ it is
$$\begin{array}{ll} \|  Q_{\cal H}(f)\|_c \,\le\,  &\kappa_{0,0} \|f\|_C +\, 2^{m-1}\, h_{x,k}\, \kappa_{1,0} \|f^{(1,0)}\|_C + \cr
\ & 2^{m-1}\, h_{y,k}\, \kappa_{0,1} \|f^{(0,1)}\|_C +\,    4^{m-1} h_{x,k} h_{y,k} \kappa_{1,1} \|f^{(1,1)}\|_C \,,  \cr
\end{array}$$
where the positive constants $\kappa_{r,s}, r, s =0,1$  are those introduced in Lemma \ref{lemmu} and  
\begin{equation} \label{Cbigdef}
C \,:=\,  \bigcup_{\ell=k_m}^k \, \bigcup_{J \in {\cal A}_\d^\ell(c)}  (\Lambda_J^\ell \cup c)\,.  
\end{equation}  
\end{lma}
\prf 
Let $c \in {\cal G}^k$. Considering the assumption on ${\cal G}_{\cal H},$  the restriction of $Q_{\cal H} (f)$ to the cell $c$ can be written as follows,
$$ Q_{\cal H} (f)_{|c}  \,=\, \sum_{\ell = k_m}^k  \sum_{J \in {\cal A}_\d^\ell(c)} \lambda_J^\ell(f) \,T_J^\ell  \,,$$
with $ {\cal A}_\d^\ell(c)$ defined as in (\ref{localex}) and where, for brevity, we have used the integer $k_m$ introduced in (\ref{kmdef}). This expression implies that $ Q_{\cal H} (f)_{|c} $ depends on $ f_{| C}\,$ where $C \supset c$ is defined as in (\ref{Cbigdef}).
Now, using the triangular inequality and considering that the THB spline basis is a nonnegative partition of unity, it can be derived that
\begin{equation} \label{locbd} 
|Q_{\cal H} (f)(x,y) |  \,\le \, \max_{k_m \le \ell \le k}\, \max_{J \in {\cal A}_\d^\ell(c)} |\lambda_J^\ell(f)| \,, \quad \forall (x,y) \in c\,.
\end{equation}

Thus the statement of the lemma is obtained by considering the functional upper bound proved in Lemma \ref{lemmu} and further uniformly upper bounding each $|\lambda_J^\ell(f)|$ in (\ref{locbd})  by replacing in (\ref{muineq}) the function norm $\| \cdot \| _{\Lambda_J^\ell}$ with $\| \cdot \| _{C}$. \QED

\noindent
The following remark outlines the relation which can be established between the size of $C$  and that of $c,$ under the assumed hypothesis on the hierarchical mesh.
\begin{rmk} \label{rmk1}
Let $R$ denote the smallest rectangle  with edges aligned with the $x$ and $y$ axes containing the convex hull of $C,$ with $C$ defined in (\ref{Cbigdef}). Further, let us denote with  $H_x(C)$ and $H_y(C)$ its $x$ and $y$ sizes, also called in the following the $x$ and $y$ diameters of $C\,.$  Then  if ${\cal G}_{\cal H}$ is admissible of class $m\,,$ there exist two positive constants $w_i > 1, i=1,2$ depending only on $\d$ and on $m$ such that for any cell $c \in {\cal G}^k, k=0,\ldots,M-1,$  it is 
\begin{equation} \label{wineq} 
H_x(C) \le\, w_1 h_{x,k}\,, \quad H_y(C) \,\le\, w_2 h_{y,k}\,.
\end{equation}
Actually it can be verified that we can   assume  
\begin{equation} \label{wdef}
w_i \,:=\, 2^{m-1}\,(d_i-1)\,, \quad i=1,2\,.
\end{equation}
\end{rmk}

Then, by using the intermediate result proved in Lemma \ref{lem2} and the inequality in (\ref{wineq}), we are ready to prove  the following main result concerning the local approximation power of $Q_{\cal H}(f) \in S_{\cal H}\,.$  
 \begin{thm} \label{locconv}
Let ${\cal G}_{\cal H}$ be admissible of class $m.$ Then there  exist three positive constants $v_i\,, i=1,2,3$   such that   $ \forall f \in C^{(d_1+1,d_2+1)}(\Omega)$ and $\forall c \in {\cal G}^k, k=0,\ldots,M-1,$ it is,
 \begin{equation} \label{thm}
\begin{array}{ll} \| f - Q_{\cal H}(f)\|_c \,\le\, &v_1\,  \|f^{(d_1+1,0)}\|_C\, h_{x,k}^{d_1+1} +   v_2\, \|f^{(0,d_2+1)}\|_C\, h_{y,k}^{d_2+1} + \cr
\ &  v_3\, \|f^{(d_1+1,d_2+1)}\|_C\,  h_{x,k}^{d_1+1} h_{y,k}^{d_2+1} \,,\cr
\end{array}
\end{equation}
where $C$ is defined in (\ref{Cbigdef}).
 \end{thm}
 
\prf
Considering that (\ref{rippol}) ensures  that $Q_{\cal H}(p) = p,\,\, \forall p \in \PP^\d,$ then the triangular inequality implies that, $\forall c \in {\cal G}^k\,,$ it is 
$$\| f - Q_{\cal H}(f)\|_c \,\le\, \|f - p\|_c \,+\, \|Q_{\cal H}(p-f)\|_c\,, \quad \forall p \in \PP^\d\,.$$
On the other hand, considering that $c \subset C$ with $C$ defined in (\ref{Cbigdef}), Lemma \ref{lem2} implies  that
$$ \begin{array}{ll} \|f - Q_{\cal H}(f) \|_c \le  & (1+\kappa_{0,0}) \|p-f\|_C +\cr
& 2^{m-1}\, h_{x,k}\, \kappa_{1,0} \|p^{(1,0)}-f^{(1,0)}\|_C + \cr \
\ & 2^{m-1}\, h_{y,k}\, \kappa_{0,1} \|p^{(0,1)}-f^{(0,1)} \|_C +\cr  
& 4^{m-1} h_{x,k} h_{y,k} \kappa_{1,1} \|p^{(1,1)}-f^{(1,1)}\|_C \,. \cr
\end{array} $$ 
Choosing $p \in \PP^\d$ as the Taylor expansion of $f$ in $\PP^\d $ at some point in $c\,,$ it is
$$\begin{array}{ll}   \| p^{(r,s)}-f^{(r,s)}\|_C \,\le\, & \|f^{(d_1+1,0)}\|_C\, H_x(C)^{d_1+1-r} +   \|f^{(0,d_2+1)}\|_C\, H_y(C)^{d_2+1-s} + \cr
\ & \|f^{(d_1+1,d_2+1)}\|_C\, H_x(C)^{d_1+1-r}  H_y(C)^{d_2+1-s}\,, \,\, r,s = 0,1\,,\cr
\end{array}
$$
where $H_x(C)$ and $H_y(C)$ denote the $x$ and $y $ diameters of $C$ introduced in Remark \ref{rmk1}.
Then the thesis is proved by using the inequality given in (\ref{wineq}) and by defining $v_i, i=1,2,3$ as follows,
$$ v_1\,:=\, K \,w_1^{d_1+1}\, , \quad v_2 \,:=\, K\, w_2^{d_2+1}\, ,\quad
 v_3 \,:=\, K \, w_1^{d_1+1}\,w_2^{d_2+1}\,,  $$
 with $K \,:=\, 1+\kappa_{0,0} + (\kappa_{1,0} +  \kappa_{0,1} + \kappa_{1,1})  \gamma\,,$ where   $\gamma := 2^{m-1} \max\{1\,,h_{x,0}\,,\,h_{y,0}\}$ and where  the constants $w_i, i=1,2\,,$ are those defined  in (\ref{wdef}). 
  \QED
  \begin{rmk} \label{meshsel}
The  statement of Theorem \ref{locconv} first of all says that $Q_{\cal H}$ mantains on a cell $c \in {\cal G}^k$  the maximal approximation order of $Q: C^{(1,1)}(\Omega) \rightarrow V^k\,.$ Furthermore,  in order to guarantee that $Q_{\cal H}$ in $S_{\cal H}$ provides the same approximation quality of $Q$ in   $V^{M-1},$  it also suggests to use larger cells  where $f^{(d_1+1,0)}$ and $f^{(0,d_2+1)}$ are almost zero and finer elsewhere.
  \end{rmk}

\section{Numerical results}\label{Tests} \label{sec_T}

\begin{figure}[H] 
%\begin{minipage}[H]{.425\textwidth}
\centering
\hspace{-.25cm}
\includegraphics[trim=0mm 0mm 0mm 0mm,clip,scale=0.30]{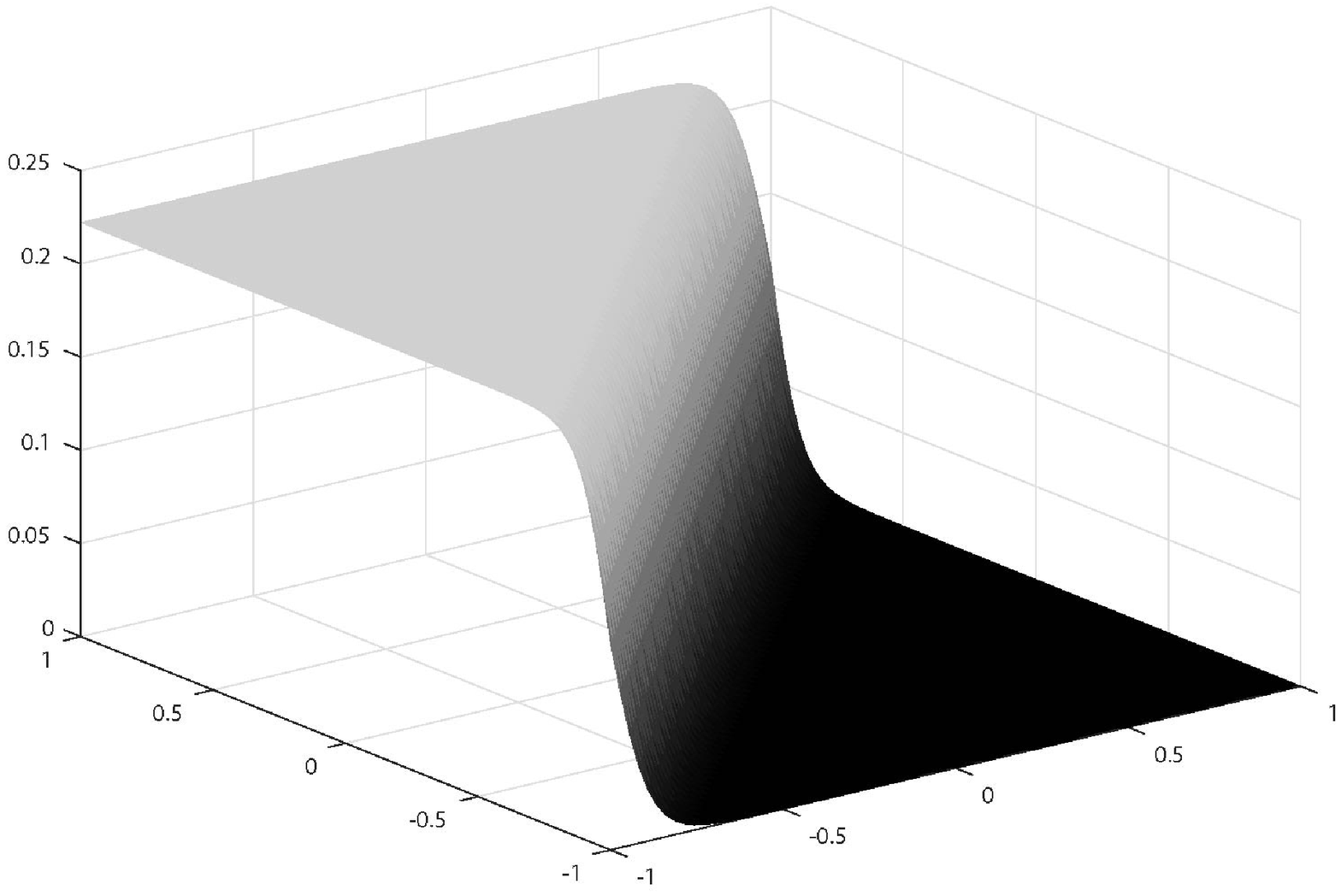}  %change figure
%\caption{Function$f_1$}
%\end{minipage} 
\hspace{-.5cm}
%\begin{minipage}[H]{.425\textwidth}
%\centering
\includegraphics[trim=0mm 0mm 0mm 0mm,clip,scale=0.40]{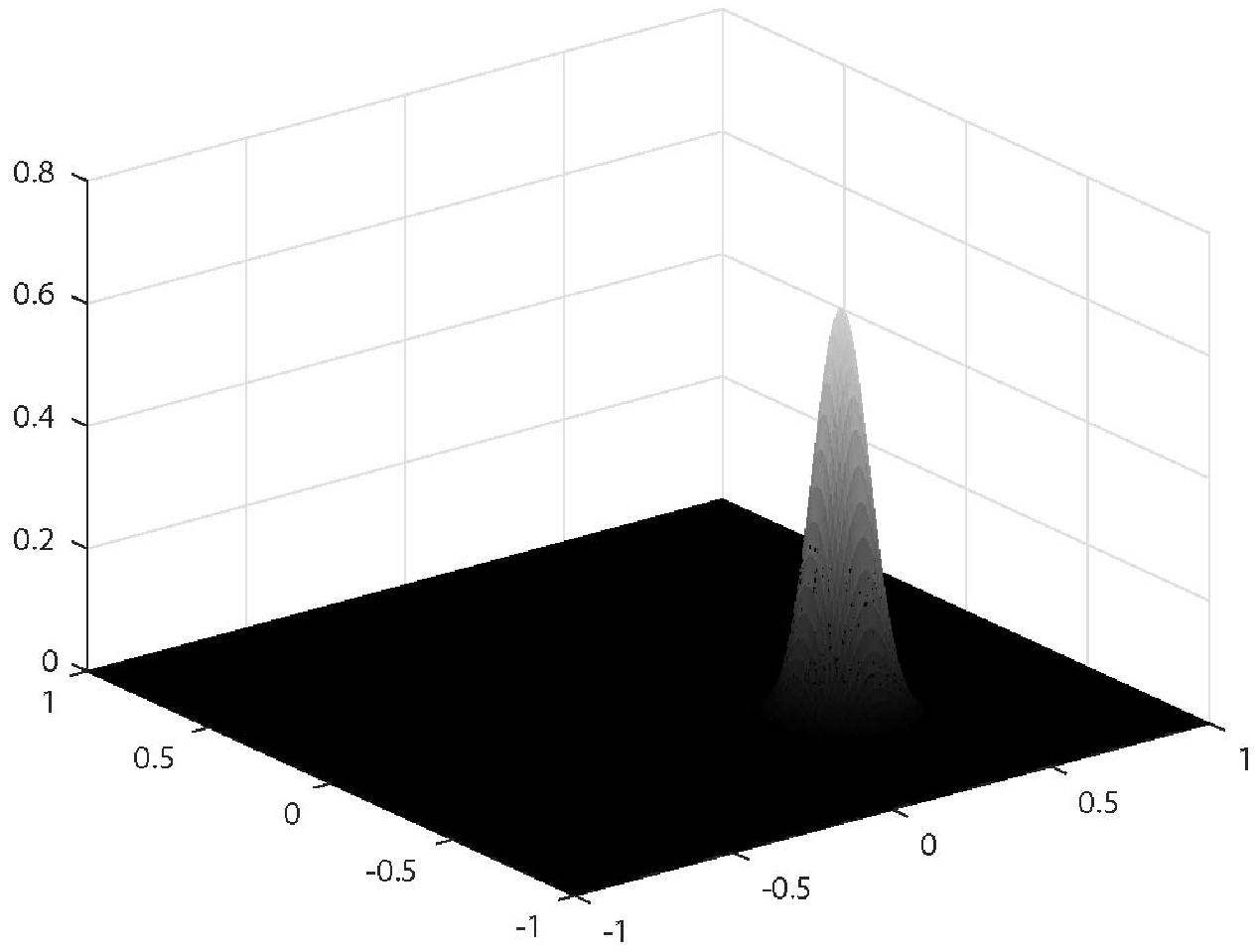}   %change figure
\caption{The test functions $f_1$ (left) and $f_2$ (right).}
%\end{minipage}  
\label{testfig}
\end{figure}

In this section we  discuss the numerical behaviour of the hierarchical BS QI operator $Q_{\cal H}$. The test functions $f(x,y)$ considered in the examples are
\begin{align*}
&f_1(x,y)=\frac{\tanh(9y-9x)+1}{9}, \qquad(x,y)\in[-1,1]^2,\\
&f_2(x,y)=\frac{2}{3\exp((10x-3)^2+(10y+4)^2)}, \qquad(x,y)\in[-1,1]^2,
\end{align*} 
see Figure \ref{testfig}.

We compare the performances of $Q_{\cal H}$ with the hierarchical QI considered in the numerical examples in \cite{SM15}. More precisely, the tensor--product version of this QI   \cite{llm2001}, here denoted by $\hat Q$, is the following,
\begin{equation*}%\label{smqi1}
\hat Q(f)(x,y):= \sum_{J \in \Gamma_{\d,\h}} \hat\lambda_J(f) \, B_J^{(\d,\h)}(x,y)\,,
\end{equation*}
where each $\hat \lambda_{J}(f)$, $J=(j,i)$, is obtained as the component $\sigma_{J}$ of the solution of the linear system
\begin{equation*}
\sum_{K:supp(B_{K}^{(\d,\h)})\cap \Theta_{J}\neq \emptyset} \sigma_{K}B_{K}^{(\d,\h)}({\bf x}^{r_1,r_2}_{J})=f({\bf x}^{r_1,r_2}_{J}), \quad r_1=0,...,d_1, \quad r_2=0,...,d_2,
\end{equation*}
with 
\begin{align*}
\Theta_{J}:=(x_{i+\lfloor d_1/2\rfloor},x_{i+\lfloor d_1/2\rfloor}+1)\times (y_{j+\lfloor d_2/2\rfloor}, y_{j+\lfloor d_2/2\rfloor+1}),\\
{\bf x}^{r_1,r_2}_{J}:=(x_{i+\lfloor d_1/2\rfloor}+\frac{r_1 h_x}{d_1},y_{j+\lfloor d_2/2\rfloor}+\frac{r_2 h_y}{d_2}). 
\end{align*}
The hierarchical extension $\hat Q_{\cal H}$ is obtained from $\hat Q$ as done for defining $Q_{\cal H}$ from $Q$ in \eqref{QH}. Essentially, the differences between $Q_{\cal H}$ and $\hat Q_{\cal H}$ lie in the type of data they use (function and derivative values for $Q_{\cal H}$ and only function values for $\hat Q_{\cal H}$), and in the number and kind of data locations.  

\smallskip
In all the tests, we considered sequences of nested hierarchical meshes on the domain $[-1,1]^2$ up to $5$ levels ($M=5$), obtained by suitably refining a uniform $8\times 8$ mesh. For all tests, the tables reporting the results obtained by applying the hierarchical QIs show the number of levels $M$, the size of the mesh $h_{M-1} :=h_{x,M-1}=h_{y,M-1}$ at the finest level, the dimension of the spline spaces and the infinity norm  $\| \cdot \|_{\Omega}$ of the errors
\begin{equation*}
e := f-Q(f),\quad e_{\cal H} := f-Q_{\cal H}(f), \quad \hat e := f-\hat Q(f), \quad \hat e_{\cal H} := f-Q_{\cal H}(f)\,.
\end{equation*}
Note that the infinity norm $\| \cdot \|_{\Omega}$ of a function has always been  numerically approximated by computing the maximum of its absolute value on a $301\times301$ uniform grid of $\Omega\,.$ Moreover, in order to have a complete understanding of the behaviour of the QIs, we also report the infinity norms of the errors with respect to the derivatives $f_x$, $f_y$ and $f_{xy}$.

\medskip
The employed hierarchical meshes are obtained by applying an automatic refinement strategy. Let $F_{\cal H}$ be the hierarchical extension, based on the mesh ${\cal G_H}$ with $M$ levels, of a tensor-product QI operator $F$, and let $P\subset \Omega$ be a set of points where the values of the approximated function $f$ (and of its derivatives, if the computation of $F$ requires them, like in the case of $Q$) are available. Let us define, for each cell $c\in {\cal G_H}$,
\begin{equation*}
\delta(F_{\cal H};c):=\max_{(x, y)\in P \cap c} \vert F_{\cal H}(x, y)-f(x, y) \vert. %[{\cal G_H}]
\end{equation*}
Then, the hierarchical mesh is generated by using the following algorithm.
\begin{itemize}
\item Start from a tensor-product mesh, which corresponds to a hierarchical mesh ${\cal G_H}$ with $1$ level ($M=1$), choose the maximum number of levels $K$ ($K=5$ in our tests) and the tolerance $\epsilon$; 
\item evaluate the corresponding $F_{\cal H}$ at the points $(x, y)\in P$; 
\item while $\delta(F_{\cal H};c)\le\epsilon$ is not satisfied for all the cells $c\in {\cal G_H}$, or $M$ is not greater than $K$, repeat the following steps:
\begin{itemize}
\item mark the cells which do not satisfy $\delta(F_{\cal H};c)\le\epsilon$;
\item for each marked cell, additionally mark the (at most) $8$ adjacent cells too;
\item obtain the new mesh ${\cal G_H}$ with $M+1$ levels by a dyadic split of the marked cells, and increase $M$ by $1$;
\item evaluate the corresponding $F_{\cal H}$ at the points $(x, y)\in P$. 
\end{itemize}
\end{itemize}
Now, in our case  the main goal is to construct a hierarchical QI operator $F_{\cal H}$ with approximation performances close to those of its tensor-product version $F$ defined using the required information on the extended lattice $\pi^e$ associated with level $K.$ Then we set $P = \pi^e \cap \Omega$ and
\begin{equation*}
\epsilon=1.5 \max_{(x, y)\in P} \vert F(x, y)-f(x, y) \vert,
\end{equation*}
which, roughly speaking, corresponds to requiring that the use of the hierarchical QI provides a maximum error at most 50\% larger than the one provided by the tensor-product version of the QI. As in our experiments $K=5$ and at the first level we use a $8 \times 8$ mesh,  $P$ is composed of the vertices of a uniform $128\times 128$ grid in $\Omega.$

Tables \ref{test1}-\ref{test1d} and \ref{test4}-\ref{test4d} report the results related to $f_1$ obtained with $(d_1,d_2)=(2,2), (3,3), (4,4)$ for $Q_{\cal H}$ and for $\hat Q_{\cal H}$, respectively.

\begin{figure} %\begin{center}
\centering
\hspace*{-0.7cm}
\includegraphics[trim=0mm 0mm 0mm 0mm,clip,scale=0.34]{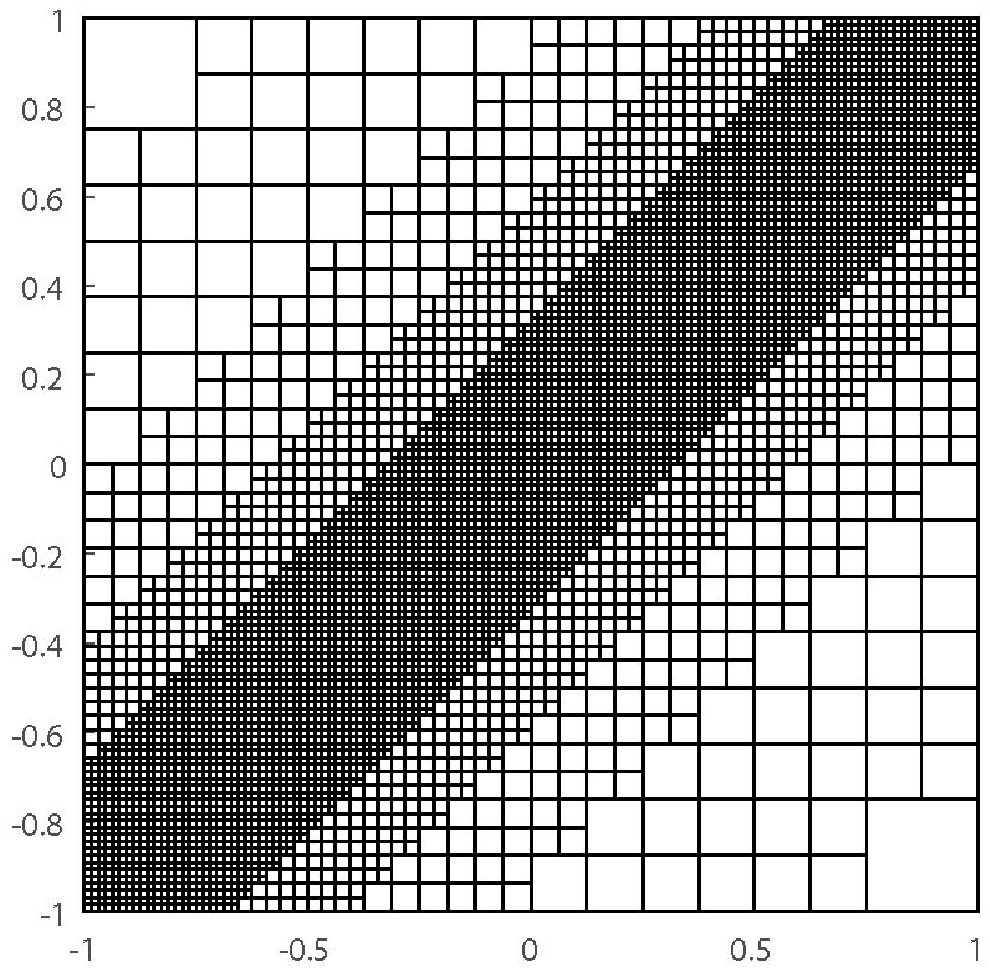}\hspace*{-0.9cm}
\includegraphics[trim=0mm 0mm 0mm 0mm,clip,scale=0.34]{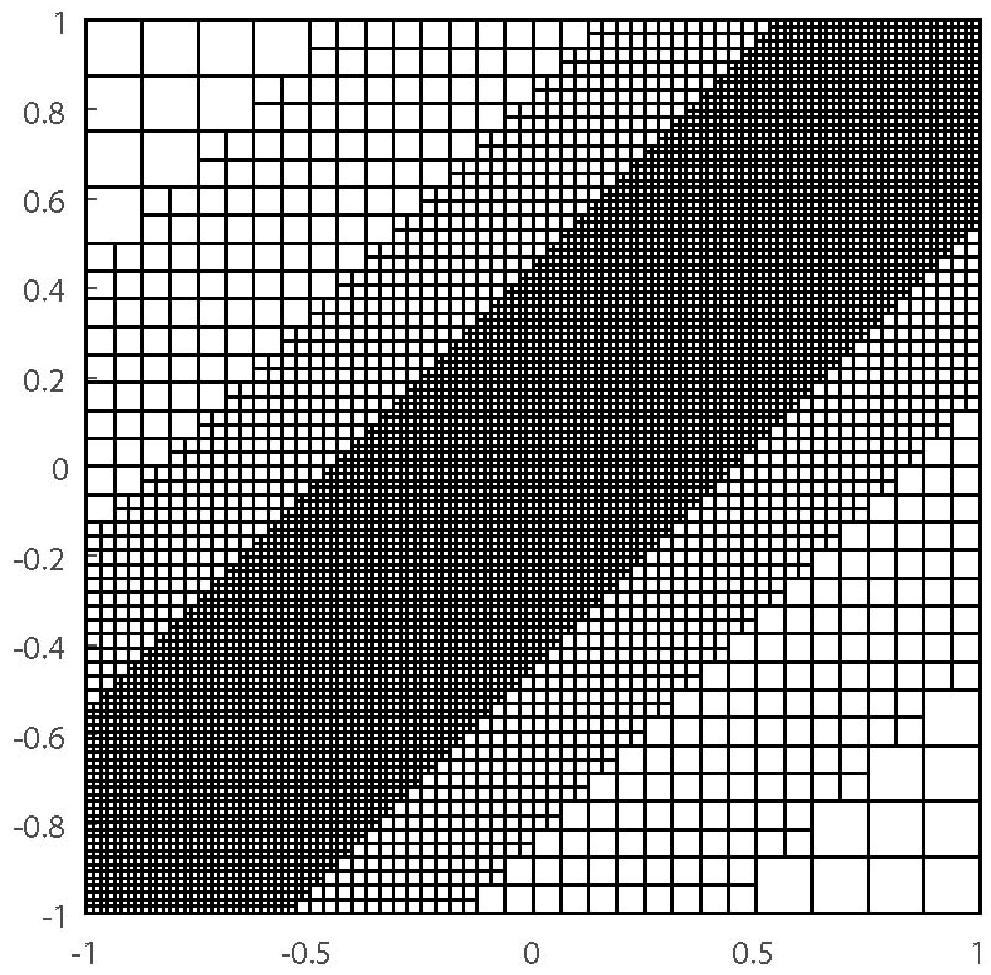}\hspace*{-0.9cm}
\includegraphics[trim=0mm 0mm 0mm 0mm,clip,scale=0.34]{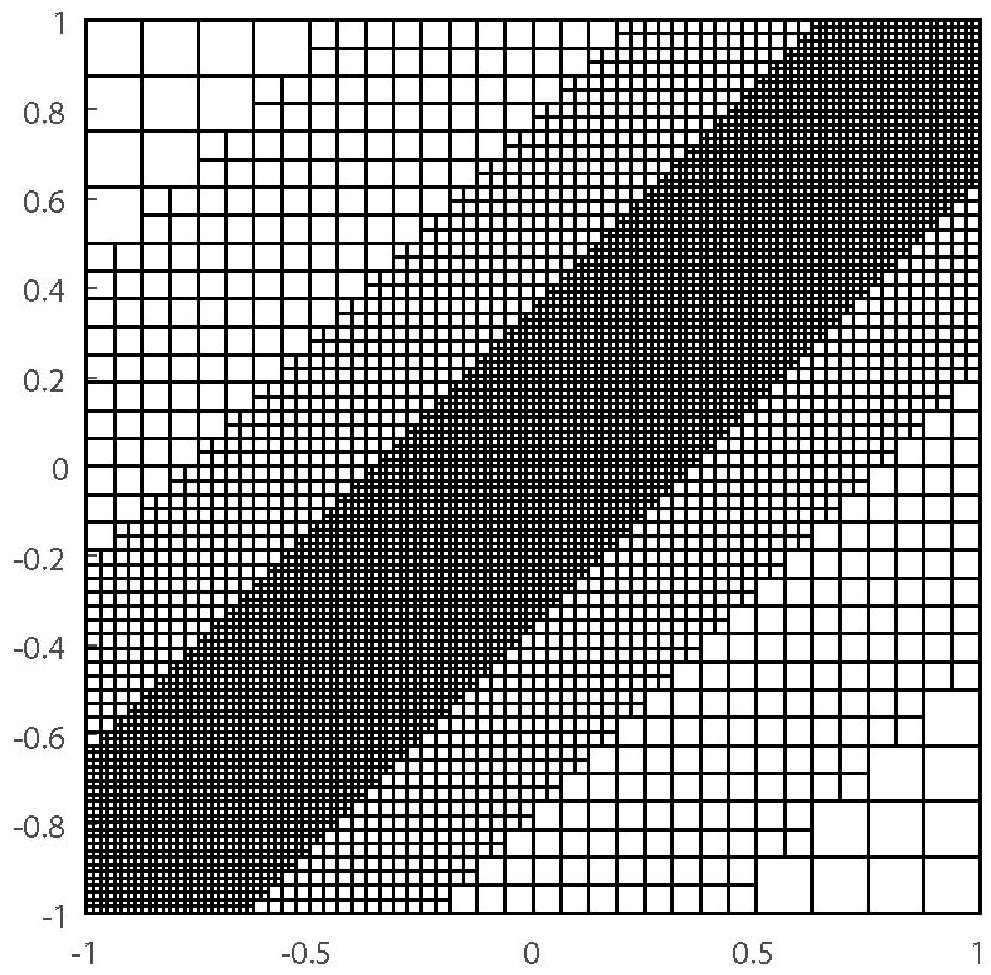}
\caption{The hierarchical meshes with 5 levels used in the approximation of $f_1$ by $Q_{\cal H}$ with (a) $(d_1,d_2)=(2,2)$, (b) $(d_1,d_2)=(3,3)$ and (c) $(d_1,d_2)=(4,4)$.}
%\end{center}
\end{figure}

\begin{figure} %\begin{center}
\centering
\hspace*{-0.7cm}
\includegraphics[trim=0mm 0mm 0mm 0mm,clip,scale=0.34]{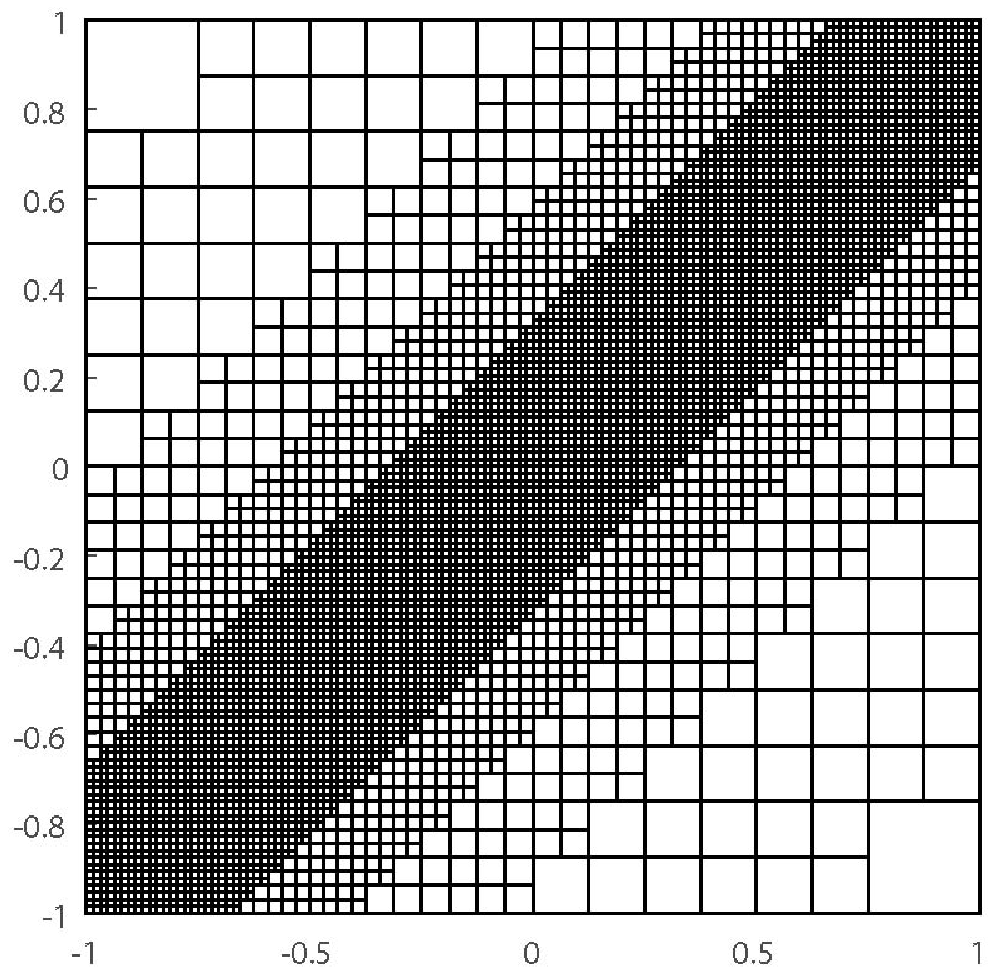}\hspace*{-0.9cm}
\includegraphics[trim=0mm 0mm 0mm 0mm,clip,scale=0.34]{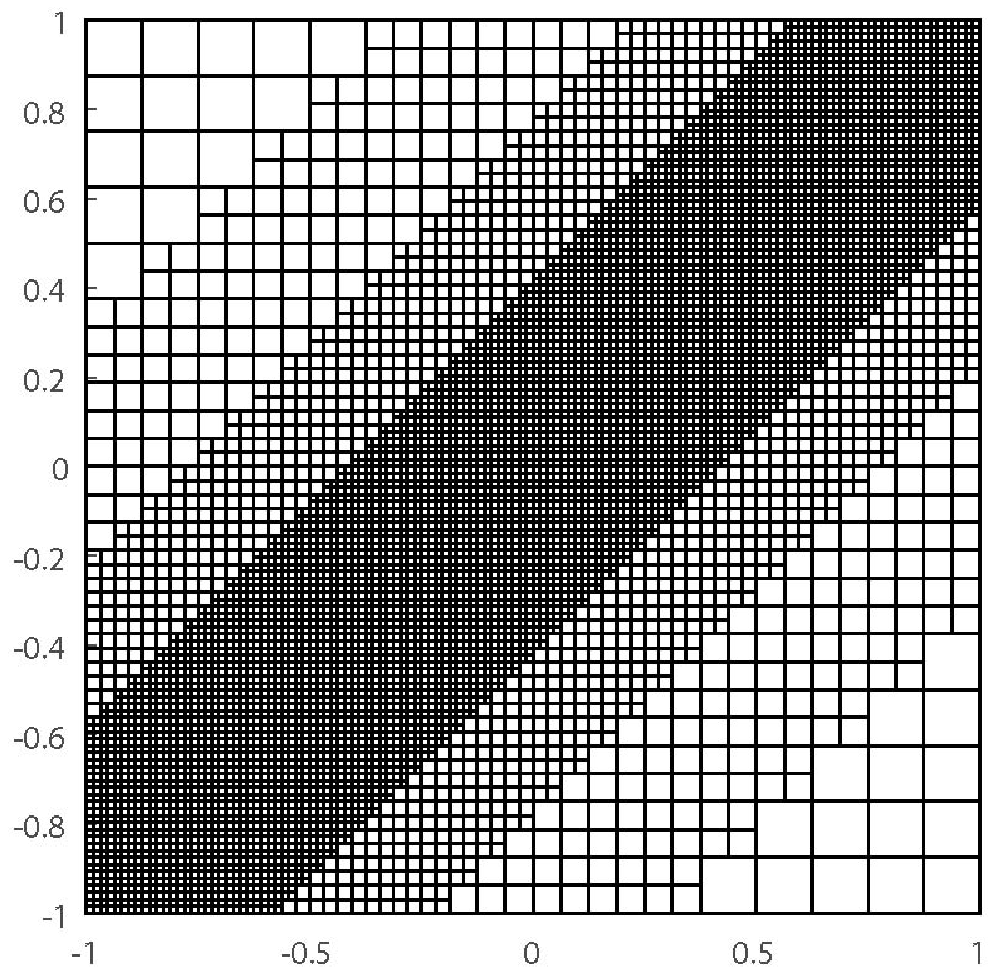}\hspace*{-0.9cm}
\includegraphics[trim=0mm 0mm 0mm 0mm,clip,scale=0.34]{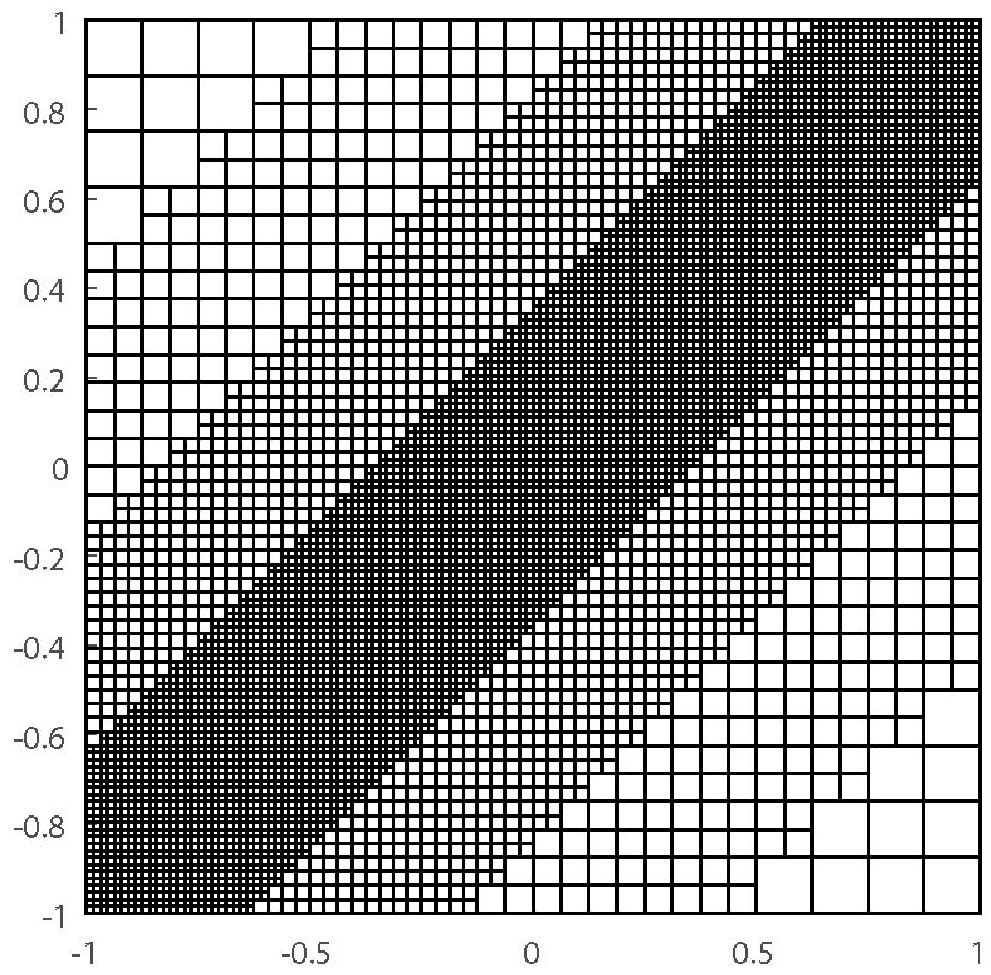}
\caption{The hierarchical meshes with 5 levels used in the approximation of $f_1$ by $\hat Q_{\cal H}$ with (a) $(d_1,d_2)=(2,2)$, (b) $(d_1,d_2)=(3,3)$ and (c) $(d_1,d_2)=(4,4)$.}
%\end{center}
\end{figure}

\begin{table}[H]
\centerline{ 
\begin{footnotesize}\begin{tabular}{SS|SS|SS} \toprule
\multicolumn{6}{c}{$(d_1,d_2)=(2,2)$}\\\midrule
    {$M$} & {${h}_{M-1} $} & {$\dim(V^{M-1}$)} & {$\Vert e\Vert_{ \Omega}$}& {$\dim(S_{\cal H})$} &  {$\Vert e_{\cal H}\Vert_{ \Omega}$} \\ \midrule
    %1  & {1/2} & {36} &  {4.796e-2}& {36} &  {4.796e-2} \\
    1  & {1/4} & {100} &  {3.050e-2} & {100} &  {3.050e-2} \\
    2  & {1/8} & {324} &  {9.982e-3} & {310} &  {9.982e-3} \\
    3  & {1/16} & {1156} &  {1.526e-3} & {862} &  {1.526e-3} \\ 
    4  & {1/32} & {4356} &  {1.312e-4} & {2368} &  {1.312e-4} \\
    5  & {1/64} & {16900} &  {1.250e-5} & {5902} &  {1.250e-5} \\ \bottomrule
\multicolumn{6}{c}{$(d_1,d_2)=(3,3)$}\\\midrule
		{$M$} & {${h}_{M-1} $} & {$\dim(V^{M-1}$)} & {$\Vert e\Vert_{ \Omega}$} & {$\dim(S_{\cal H})$} &  {$\Vert e_{\cal H}\Vert_{ \Omega}$} \\ \midrule
   % 1  & {1/2} & {49}  & {1.119e-1} & {49}  & {1.119e-1} \\
    1  & {1/4} & {121}  & {4.581e-2} & {121}  & {4.581e-2} \\
    2  & {1/8} & {361} & {8.168e-3} & {361}  & {8.168e-3} \\
    3  & {1/16} & {1225}  & {5.951e-4} & {1117}  & {5.951e-4} \\ 
    4  & {1/32} & {4489}  & {2.414e-5} & {3139}  & {2.414e-5} \\
    5  & {1/64} & {17161}  & {1.115e-6} & {7873}  & {1.115e-6} \\ \bottomrule
\multicolumn{6}{c}{$(d_1,d_2)=(4,4)$}\\\midrule
{$M$} & {${h}_{M-1} $} & {$\dim(V^{M-1}$)} & {$\Vert e\Vert_{ \Omega}$} & {$\dim(S_{\cal H})$} & {$\Vert e_{\cal H}\Vert_{ \Omega}$} \\ \midrule
    %1  & {1/2} & {64}  & {1.207e-1} & {64} & {1.207e-1} \\
    1  & {1/4} & {144}  & {6.842e-2} & {144} & {6.842e-2} \\
    2  & {1/8} & {400}  & {1.034e-2} & {400} & {1.034e-2} \\
    3  & {1/16} & {1296}  & {3.980e-4} & {1172} & {3.980e-4} \\ 
    4  & {1/32} & {4624}  & {8.828e-6} & {3056} & {8.828e-6} \\
    5  & {1/64} & {17424}  & {1.512e-7} & {6756} & {1.512e-7} \\ \bottomrule
\end{tabular}\end{footnotesize}
}
\caption{Numerical results for $Q_{\cal H}$ tested on $f_1$ and compared with $Q$.}\label{test1}
\end{table}

\begin{table}[H]%\begin{tiny}
\centerline{ 
\begin{footnotesize}\begin{tabular}{S|SSS|SSS} \toprule
\multicolumn{7}{c}{$(d_1,d_2)=(2,2)$}\\\midrule
    {$M$} & {$\Vert e_x\Vert_{ \Omega}$} & {$\Vert e_y\Vert_{ \Omega}$} & {$\Vert e_{xy}\Vert_{ \Omega}$}& {$\Vert (e_{\cal H})_x\Vert_{ \Omega}$} & {$\Vert (e_{\cal H})_y\Vert_{ \Omega}$} & {$\Vert (e_{\cal H})_{xy}\Vert_{ \Omega}$} \\ \midrule
    %1  & {6.593e-1} & {6.593e-1} &  {6.764e-0}& {6.593e-1} & {6.593e-1} &  {6.764e-0} \\
    1  & {4.933e-1} & {4.933e-1} &  {6.185e-0}& {4.933e-1} & {4.933e-1} &  {6.185e-0} \\
    2  & {2.218e-1} & {2.218e-1} &  {4.133e-0} & {2.218e-1} & {2.218e-1} &  {4.133e-0} \\
    3  & {5.266e-2} & {5.266e-2} &  {1.537e-0} & {5.266e-2} & {5.266e-2} &  {1.537e-0} \\ 
    4  & {1.017e-2} & {1.017e-2} &  {3.019e-1} & {1.017e-2} & {1.017e-2} &  {3.019e-1} \\
    5  & {3.088e-3} & {3.088e-3} &  {1.113e-1} & {3.088e-3} & {3.088e-3} &  {1.113e-1} \\ \bottomrule
		\multicolumn{7}{c}{$(d_1,d_2)=(3,3)$}\\\midrule
		{$M$} & {$\Vert e_x\Vert_{ \Omega}$} & {$\Vert e_y\Vert_{ \Omega}$} & {$\Vert e_{xy}\Vert_{ \Omega}$}& {$\Vert (e_{\cal H})_x\Vert_{ \Omega}$} & {$\Vert (e_{\cal H})_y\Vert_{ \Omega}$} & {$\Vert (e_{\cal H})_{xy}\Vert_{ \Omega}$} \\ \midrule
    %1  & {1.027e-0} & {1.027e-0} &  {7.081e-0}& {1.027e-0} & {1.027e-0} &  {7.081e-0} \\
    1  & {6.339e-1} & {6.339e-1} &  {6.600e-0}& {6.339e-1} & {6.339e-1} &  {6.600e-0} \\
    2  & {1.812e-1} & {1.812e-1} &  {3.741e-0} & {1.812e-1} & {1.812e-1} &  {3.741e-0} \\
    3  & {1.835e-2} & {1.835e-2} &  {7.533e-1} & {1.835e-2} & {1.835e-2} &  {7.533e-1} \\ 
    4  & {1.263e-3} & {1.263e-3} &  {7.065e-2} & {1.263e-3} & {1.263e-3} &  {7.065e-2} \\
    5  & {9.971e-5} & {9.971e-5} &  {6.179e-3} & {9.972e-5} & {9.972e-5} &  {6.179e-3} \\ \bottomrule
		\multicolumn{7}{c}{$(d_1,d_2)=(4,4)$}\\\midrule
		{$M$} & {$\Vert e_x\Vert_{ \Omega}$} & {$\Vert e_y\Vert_{ \Omega}$} & {$\Vert e_{xy}\Vert_{ \Omega}$}& {$\Vert (e_{\cal H})_x\Vert_{ \Omega}$} & {$\Vert (e_{\cal H})_y\Vert_{ \Omega}$} & {$\Vert (e_{\cal H})_{xy}\Vert_{ \Omega}$} \\ \midrule
    %1  & {1.084e-0} & {1.084e-0} &  {7.200e-0}& {1.084e-0} & {1.084e-0} &  {7.200e-0} \\
    1  & {8.318e-1} & {8.318e-1} &  {7.401e-0}& {8.318e-1} & {8.318e-1} &  {7.401e-0} \\
    2  & {2.212e-1} & {2.212e-1} &  {4.012e-0} & {2.212e-1} & {2.212e-1} &  {4.012e-0} \\
    3  & {1.457e-2} & {1.457e-2} &  {5.285e-1} & {1.457e-2} & {1.457e-2} &  {5.285e-1} \\ 
    4  & {4.846e-4} & {4.846e-4} &  {2.389e-2} & {4.846e-4} & {4.846e-4} &  {2.389e-2} \\
    5  & {1.401e-5} & {1.401e-5} &  {6.941e-4} & {1.401e-5} & {1.401e-5} &  {6.941e-4} \\ \bottomrule
\end{tabular}\end{footnotesize}
}%\end{tiny}
\caption{Errors with respect to the derivatives for $Q_{\cal H}$ tested on $f_1$ and compared with $Q$.}\label{test1d}
\end{table}

\begin{table}[H]
\centerline{ 
\begin{footnotesize}\begin{tabular}{SS|SS|SS} \toprule
\multicolumn{6}{c}{$(d_1,d_2)=(2,2)$}\\\midrule
     {$M$} & {${h}_{M-1} $} & {$\dim(V^{M-1}$)} & {$\Vert \hat e\Vert_{ \Omega}$}& {$\dim(S_{\cal H})$} &  {$\Vert \hat e_{\cal H}\Vert_{ \Omega}$} \\ \midrule
    %1  & {1/2} & {36} & {7.433e-2} & {36} & {7.433e-2} \\
    1  & {1/4} & {100} & {4.165e-2} & {100} & {4.165e-2} \\
    2  & {1/8} & {324} & {1.007e-2} & {310} & {1.007e-2} \\
    3  & {1/16} & {1156} & {1.292e-3} & {862} & {1.292e-3} \\ 
    4  & {1/32} & {4356} & {1.160e-4} & {2368} & {1.160e-4} \\
    5  & {1/64} & {16900} & {1.109e-5} & {5902} & {1.109e-5} \\ \bottomrule
		\multicolumn{6}{c}{$(d_1,d_2)=(3,3)$}\\\midrule
		{$M$} & {${h}_{M-1} $} & {$\dim(V^{M-1}$)} & {$\Vert \hat e\Vert_{ \Omega}$}& {$\dim(S_{\cal H})$} &  {$\Vert \hat e_{\cal H}\Vert_{ \Omega}$} \\ \midrule
    %1  & {1/2} & {49}  & {4.699e-2} & {49} & {4.699e-2} \\
    1  & {1/4} & {121} & {7.646e-2} & {121} & {7.646e-2} \\
    2  & {1/8} & {361} & {1.290e-2} & {361} & {1.290e-2} \\
    3  & {1/16} & {1225} & {4.584e-4} & {1075} & {4.584e-4} \\ 
    4  & {1/32} & {4489} & {1.123e-5} & {2971}  & {1.123e-5} \\
    5  & {1/64} & {17161} & {9.183e-7} & {7393} & {9.183e-7} \\ \bottomrule
		\multicolumn{6}{c}{$(d_1,d_2)=(4,4)$}\\\midrule
		{$M$} & {${h}_{M-1} $} & {$\dim(V^{M-1}$)} & {$\Vert \hat e\Vert_{ \Omega}$}& {$\dim(S_{\cal H})$} &  {$\Vert \hat e_{\cal H}\Vert_{ \Omega}$} \\ \midrule
    %1  & {1/2} & {64} & {1.730e-1} & {64} & {1.730e-1} \\
    1  & {1/4} & {144} & {1.015e-1} & {144}  & {1.015e-1} \\
    2  & {1/8} & {400} & {2.945e-2} & {400} & {2.945e-2} \\
    3  & {1/16} & {1296} & {4.799e-4} & {1172} & {4.799e-4} \\ 
    4  & {1/32} & {4624} & {1.061e-5} & {3186} & {1.061e-5} \\
    5  & {1/64} & {17424} & {1.980e-7} & {6886} & {1.980e-7} \\ \bottomrule
\end{tabular}\end{footnotesize}
}
\caption{Numerical results for $\hat Q_{\cal H}$ tested on $f_1$ and compared with $\hat Q$.}\label{test4}
\end{table}

\begin{table}[H]
\centerline{ 
\begin{footnotesize}\begin{tabular}{S|SSS|SSS} \toprule
\multicolumn{7}{c}{$(d_1,d_2)=(2,2)$}\\\midrule
    {$M$} & {$\Vert \hat e_x\Vert_{ \Omega}$} & {$\Vert \hat e_y\Vert_{ \Omega}$} & {$\Vert \hat e_{xy}\Vert_{ \Omega}$}& {$\Vert (\hat e_{\cal H})_x\Vert_{ \Omega}$} & {$\Vert (\hat e_{\cal H})_y\Vert_{ \Omega}$} & {$\Vert (\hat e_{\cal H})_{xy}\Vert_{ \Omega}$} \\ \midrule
    %1  & {8.401e-1} & {8.401e-1} &  {6.892e-0} & {8.401e-1} & {8.401e-1} &  {6.892e-0} \\
    1  & {5.988e-1} & {5.988e-1} &  {6.451e-0} & {5.988e-1} & {5.988e-1} &  {6.451e-0} \\
    2  & {2.229e-1} & {2.229e-1} &  {4.127e-0} & {2.229e-1} & {2.229e-1} &  {4.127e-0} \\
    3  & {4.645e-2} & {4.645e-2} &  {1.425e-0} & {4.645e-2} & {4.645e-2} &  {1.425e-0} \\ 
    4  & {1.101e-2} & {1.101e-2} &  {3.193e-1} & {1.101e-2} & {1.101e-2} &  {3.193e-1} \\
    5  & {3.152e-3} & {3.152e-3} &  {1.150e-1} & {3.152e-3} & {3.152e-3} &  {1.150e-1} \\ \bottomrule
		\multicolumn{7}{c}{$(d_1,d_2)=(3,3)$}\\\midrule
		{$M$} & {$\Vert \hat e_x\Vert_{ \Omega}$} & {$\Vert \hat e_y\Vert_{ \Omega}$} & {$\Vert \hat e_{xy}\Vert_{ \Omega}$}& {$\Vert (\hat e_{\cal H})_x\Vert_{ \Omega}$} & {$\Vert (\hat e_{\cal H})_y\Vert_{ \Omega}$} & {$\Vert (\hat e_{\cal H})_{xy}\Vert_{ \Omega}$} \\ \midrule
    %1  & {6.500e-0} & {6.500e-0} &  {6.781e-0}& {6.500e-0} & {6.500e-0} &  {6.781e-0} \\
    1  & {8.666e-1} & {8.666e-1} &  {7.183e-0}& {8.666e-1} & {8.666e-1} &  {7.183e-0} \\
    2  & {2.525e-1} & {2.525e-1} &  {4.271e-0} & {2.525e-1} & {2.525e-1} &  {4.271e-0} \\
    3  & {1.709e-2} & {1.709e-2} &  {6.800e-1} & {1.709e-2} & {1.709e-2} &  {6.800e-1} \\ 
    4  & {9.541e-4} & {9.541e-4} &  {6.068e-2} & {9.541e-4} & {9.541e-4} &  {6.068e-2} \\
    5  & {1.001e-4} & {1.001e-4} &  {6.504e-3} & {1.001e-4} & {1.001e-4} &  {6.504e-3} \\ \bottomrule
		\multicolumn{7}{c}{$(d_1,d_2)=(4,4)$}\\\midrule
		{$M$} & {$\Vert \hat e_x\Vert_{ \Omega}$} & {$\Vert \hat e_y\Vert_{ \Omega}$} & {$\Vert \hat e_{xy}\Vert_{ \Omega}$}& {$\Vert (\hat e_{\cal H})_x\Vert_{ \Omega}$} & {$\Vert (\hat e_{\cal H})_y\Vert_{ \Omega}$} & {$\Vert (\hat e_{\cal H})_{xy}\Vert_{ \Omega}$} \\ \midrule
    %1  & {5.397e-1} & {5.397e-1} &  {6.553e-0}& {5.397e-1} & {5.397e-1} &  {6.553e-0} \\
    1  & {5.247e-1} & {5.247e-1} &  {5.096e-0}& {5.247e-1} & {5.247e-1} &  {5.096e-0} \\
    2  & {4.704e-1} & {4.704e-1} &  {5.823e-0} & {4.704e-1} & {4.704e-1} &  {5.823e-0} \\
    3  & {1.687e-2} & {1.687e-2} &  {5.774e-1} & {1.687e-2} & {1.687e-2} &  {5.774e-1} \\ 
    4  & {5.538e-4} & {5.538e-4} &  {2.597e-2} & {5.538e-4} & {5.538e-4} &  {2.597e-2} \\
    5  & {1.623e-5} & {1.623e-5} &  {7.786e-4} & {1.623e-5} & {1.623e-5} &  {7.786e-4} \\ \bottomrule
\end{tabular}\end{footnotesize}
}
\caption{Errors with respect to the derivatives for $\hat Q_{\cal H}$ tested on $f_1$ and compared with $\hat Q$.}\label{test4d}
\end{table}
We further compared $Q_{\cal H}$ and $\hat Q_{\cal H}$ by examining the functional evaluations needed in the two cases: as shown in Table \ref{eval1}, while for $(d_1,d_2)=(2,2)$ the number of required evaluations  is smaller for $\hat Q_{\cal H}$, when $(d_1,d_2)=(3,3),(4,4)$ $\hat Q_{\cal H}$ requires significantly more evaluations than $Q_{\cal H}$, in spite of using, in some cases, slightly less refined hierarchical meshes (which lead to spline spaces with smaller dimension). Note that for the corresponding tensor--product versions, the number of evaluations is  $4(N_1+2d_1-1)(N_2+2d_2-1)$ for $Q$ and $(d_1(N_1+d_1)+1)(d_2(N_2+d_2)+1)$ for $\hat Q$, which is the reason behind the increasing difference between the evaluations needed by $Q_{\cal H}$ and $\hat Q_{\cal H}$ as $d_1$ and $d_2$ grow.

Tables \ref{test8}-\ref{test8d} and \ref{test11}-\ref{test11d} report, respectively, the results for $Q_{\cal H}$ and for $\hat Q_{\cal H}$ applied to $f_2$ with $(d_1,d_2)=(3,3)$. In this case, we can observe that at the coarsest levels $Q_{\cal H}$ and $\hat Q_{\cal H}$ behave differently, with $Q_{\cal H}$ showing a better accuracy. Note that such differences are not related to the hierarchical nature of the QIs, since the same behaviour can be observed in their tensor--product versions (as in the previous case, the choice of the hierarchical meshes allows to obtain essentially the same maximum error as in the tensor--product case).

\begin{table}
\centerline{ 
\begin{footnotesize}\begin{tabular}{S|SSSSS} \toprule
    \multicolumn{6}{c}{$(d_1,d_2)=(2,2)$}\\\midrule
		{$M$}&{1} & {2} & {3} & {4} & {5} \\
    {$Q$} & {484} & {1444} & {4900} & {17956} & {68644} \\
		{$\hat Q$} & {441} & {1369} & {4761} & {17689} & {68121}\\
		{$Q_{\cal H}$} & {484} & {1404} & {3756} & {10052} & {24716}\\
		{$\hat Q_{\cal H}$} & {441} & {1299} & {3501} & {9471} & {23433}\\\bottomrule	
		\multicolumn{6}{c}{$(d_1,d_2)=(3,3)$}\\\midrule %{70756} & {155236} & {21348} & {41508}
		{$M$}&{1} & {2} & {3} & {4} & {5} \\
    {$Q$} & {676} & {1764} & {5476} & {19044} & {70756} \\
		{$\hat Q$} & {1156} & {3364} & {11236} & {40804} & {155236}\\
		{$Q_{\cal H}$} & {676} & {1764} & {5076} & {13708} & {33700}\\
		{$\hat Q_{\cal H}$} & {1156} & {3364} & {9732} & {26498} & {65446}\\\bottomrule	
		\multicolumn{6}{c}{$(d_1,d_2)=(4,4)$}\\\midrule %{72900} & {279841} & {24108} & {77713}
		{$M$}&{1} & {2} & {3} & {4} & {5}\\
    {$Q$} & {900} & {2116} & {6084} & {20166} & {72900} \\
		{$\hat Q$} & {2401} & {6561} & {21025} & {74529} & {279841}\\
		{$Q_{\cal H}$} & {900} & {2116} & {5684} & {14084} & {30516}\\
		{$\hat Q_{\cal H}$} & {2401} & {6561} & {18641} & {49825} & {106065}\\\bottomrule	
\end{tabular}\end{footnotesize}
}
\caption{Number of functional evaluations needed by $Q$, $\hat Q$, $Q_{\cal H}$ and $\hat Q_{\cal H}$ in the test on $f_1$.}\label{eval1}
\end{table}

\begin{figure} %\begin{center}
%\hspace*{-0.7cm}
\centering
\includegraphics[trim=0mm 0mm 0mm 0mm,clip,scale=0.41]{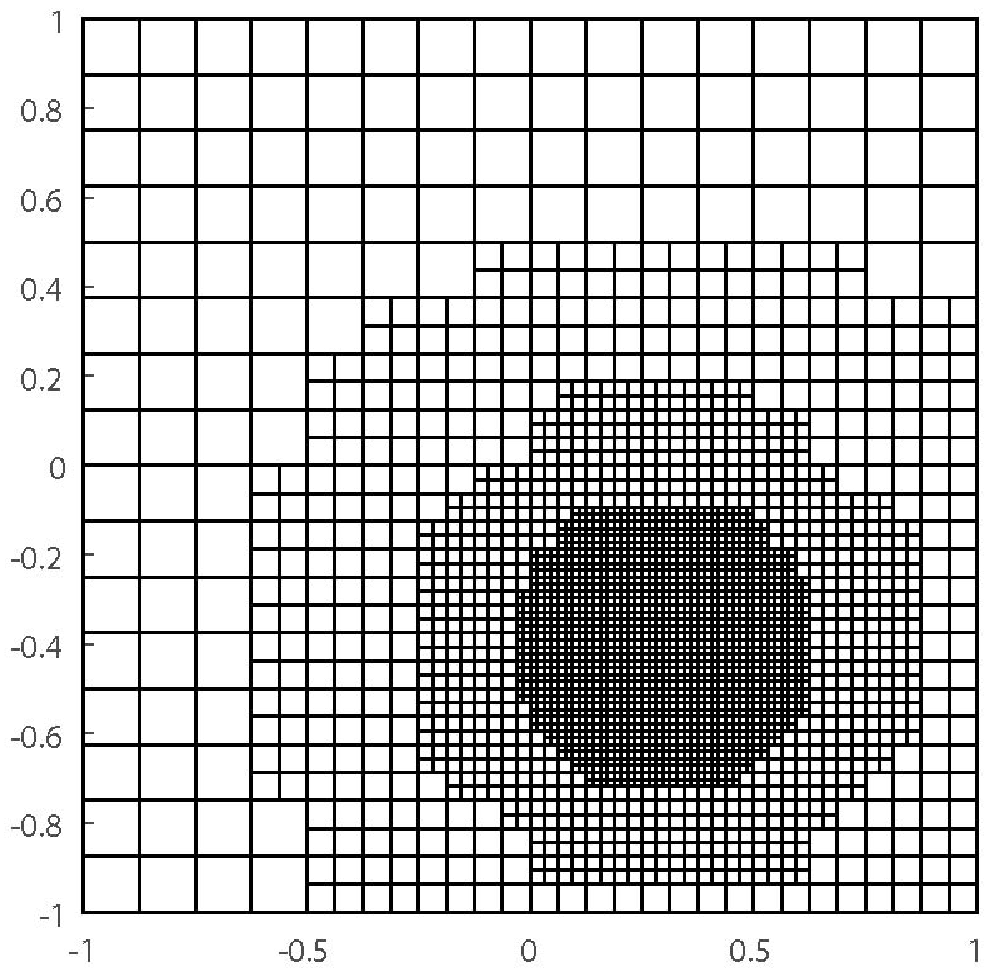}%\hspace*{-0.9cm}
\includegraphics[trim=0mm 0mm 0mm 0mm,clip,scale=0.41]{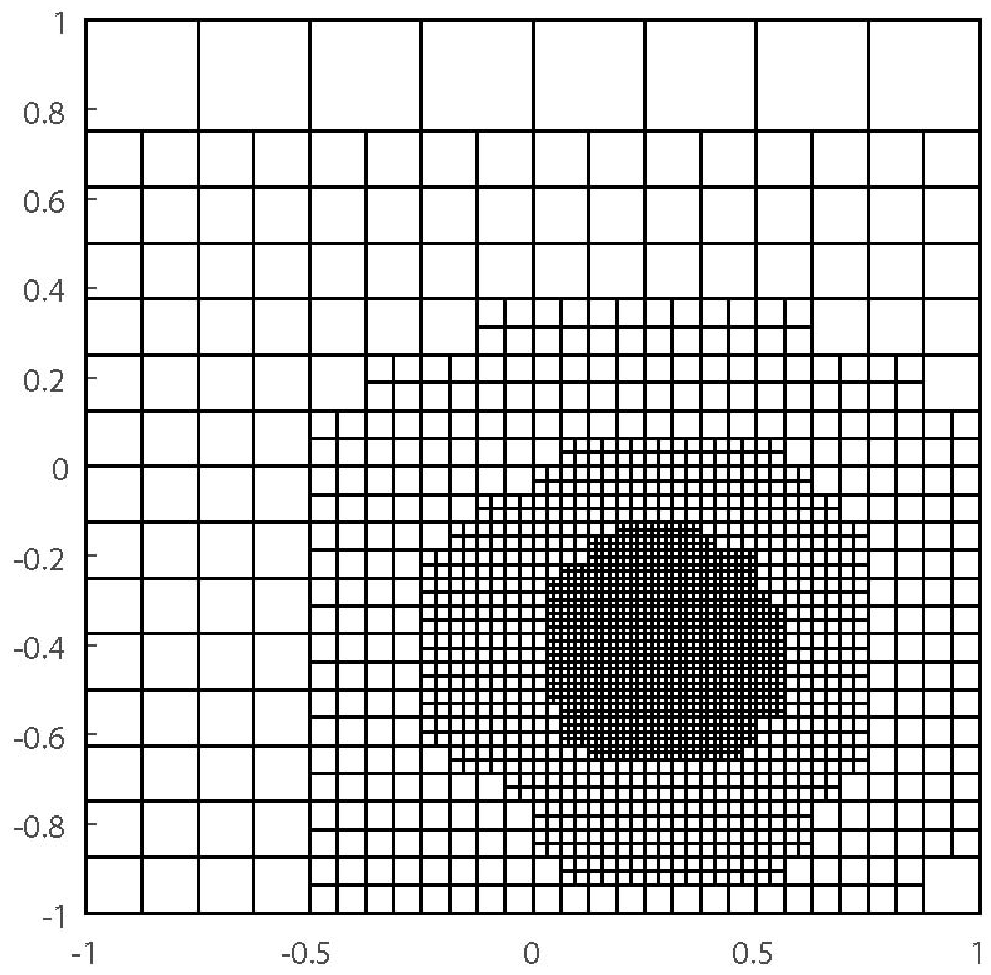}
\caption{The hierarchical meshes with 5 levels used to approximate $f_2$ by $Q_{\cal H}$ (right) and by $\hat Q_{\cal H}$ (left) with $(d_1,d_2)=(3,3)$.}
%\end{center}
\end{figure}

\begin{table}[H]
\centerline{ 
\begin{footnotesize}\begin{tabular}{SS|SS|SS} \toprule
    {$M$} & {${h}_{M-1} $} & {$\dim(V^{M-1}$)} & {$\Vert e\Vert_{ \Omega}$} & {$\dim(S_{\cal H})$} & {$\Vert e_{\cal H}\Vert_{ \Omega}$} \\ \midrule
    %1  & {1/2} & {49} & {6.681e-1} & {49} & {6.681e-1} \\
    1  & {1/4} & {121} & {5.763e-1} & {121} & {5.763e-1} \\
    2  & {1/8} & {361} & {1.974e-1} & {361} & {1.974e-1} \\
    3  & {1/16} & {1225} & {1.662e-2} & {787} & {1.662e-2} \\ 
    4  & {1/32} & {4489} & {6.559e-4} & {1471} & {6.559e-4} \\
    5  & {1/64} & {17161} & {2.760e-5} & {2440} & {2.760e-5} \\ \bottomrule
\end{tabular}\end{footnotesize}
}
\caption{Numerical results for $Q_{\cal H}$ tested on $f_2$ with $(d_1,d_2)=(3,3)$ and compared with $Q$.}\label{test8}
\end{table}

\begin{table}[H]
\centerline{ 
\begin{footnotesize}\begin{tabular}{S|SSS|SSS} \toprule
    {$M$} & {$\Vert e_x\Vert_{ \Omega}$} & {$\Vert e_y\Vert_{ \Omega}$} & {$\Vert e_{xy}\Vert_{ \Omega}$}& {$\Vert (e_{\cal H})_x\Vert_{ \Omega}$} & {$\Vert (e_{\cal H})_y\Vert_{ \Omega}$} & {$\Vert (e_{\cal H})_{xy}\Vert_{ \Omega}$} \\ \midrule
    %1  & {5.738e-0} & {5.721e-0} &  {4.912e+1} & {5.738e-0} & {5.721e-0} &  {4.912e+1} \\
    1  & {5.732e-0} & {6.403e-0} &  {5.385e+1} & {5.732e-0} & {6.403e-0} &  {5.385e+1} \\
    2  & {3.504e-0} & {2.585e-0} &  {3.181e+1} & {3.504e-0} & {2.585e-0} &  {3.181e+1} \\
    3  & {4.127e-1} & {4.067e-1} &  {4.762e-0} & {4.127e-1} & {4.067e-1} &  {4.762e-0} \\ 
    4  & {2.581e-2} & {2.620e-2} &  {2.736e-1} & {2.581e-2} & {2.620e-2} &  {2.736e-1} \\
    5  & {2.531e-3} & {2.537e-3} &  {2.414e-2} & {2.531e-3} & {2.537e-3} &  {2.414e-2} \\ \bottomrule
\end{tabular}\end{footnotesize}
}
\caption{Errors with respect to the derivatives for $Q_{\cal H}$ tested on $f_2$ with $(d_1,d_2)=(3,3)$ and compared with $Q$.}\label{test8d}
\end{table}

\begin{table}[H]
\centerline{ 
\begin{footnotesize}\begin{tabular}{SS|SS|SS} \toprule
    {$M$} & {${h}_{M-1} $} & {$\dim(V^{M-1}$)} & {$\Vert \hat e\Vert_{ \Omega}$} & {$\dim(S_{\cal H})$} & {$\Vert \hat e_{\cal H}\Vert_{ \Omega}$} \\ \midrule
    %1  & {1/2} & {49} & {5.071e-0} & {49} & {5.071e-0} \\
    1  & {1/4} & {121} & {3.087e-0} & {121} & {3.087e-0} \\
    2  & {1/8} & {361} & {2.285e-1} & {310} & {2.285e-1} \\
    3  & {1/16} & {1225} & {1.311e-2} & {667} & {1.311e-2} \\ 
    4  & {1/32} & {4489} & {6.365e-4} & {1210} & {6.365e-4} \\
    5 & {1/64} & {17161} & {3.154e-5} & {1846} & {3.154e-5} \\ \bottomrule
\end{tabular}\end{footnotesize}
}
\caption{Numerical results for $\hat Q_{\cal H}$ tested on $f_2$ with $(d_1,d_2)=(3,3)$ and compared with $\hat Q$.}\label{test11}
\end{table}

\begin{table}[H]
\centerline{ 
\begin{footnotesize}\begin{tabular}{S|SSS|SSS} \toprule
    {$M$} & {$\Vert \hat e_x\Vert_{ \Omega}$} & {$\Vert \hat e_y\Vert_{ \Omega}$} & {$\Vert \hat e_{xy}\Vert_{ \Omega}$}& {$\Vert (\hat e_{\cal H})_x\Vert_{ \Omega}$} & {$\Vert (\hat e_{\cal H})_y\Vert_{ \Omega}$} & {$\Vert (\hat e_{\cal H})_{xy}\Vert_{ \Omega}$} \\ \midrule
    %1  & {1.501e+1} & {1.084e+1} &  {6.449e+1} & {1.501e+1} & {1.084e+1} &  {6.449e+1} \\
    1  & {1.203e+1} & {1.300e+1} &  {7.423e+1} & {1.203e+1} & {1.300e+1} &  {7.423e+1} \\
    2  & {4.189e-0} & {2.944e-0} &  {3.652e+1} & {4.189e-0} & {2.944e-0} &  {3.652e+1} \\
    3  & {4.064e-1} & {4.744e-1} &  {5.057e-0} & {4.064e-1} & {4.744e-1} &  {5.057e-0} \\ 
    4  & {3.313e-2} & {3.374e-2} &  {3.442e-1} & {3.313e-2} & {3.374e-2} &  {3.442e-1} \\
    5  & {2.961e-3} & {2.714e-3} &  {2.918e-2} & {2.961e-3} & {2.714e-3} &  {5.931e-2} \\ \bottomrule
\end{tabular}\end{footnotesize}
}
\caption{Errors with respect to the derivatives for $\hat Q_{\cal H}$ tested on $f_2$ with $(d_1,d_2)=(3,3)$ and compared with $\hat Q$.}\label{test11d}
\end{table}

Finally, we present two test where we approximated $f_1$ and $f_2$ with a modified version of the operator $Q_{\cal H}$, denoted by $Q^a_{\cal H}$, which does not require the derivative values. In this case the computation of  each functional coefficient $\lambda_J^{\ell}(f)$ defining in (\ref{QH}) the hierarchical quasi--interpolant is done by replacing the derivative values  with their suitable finite-differences approximations which are locally defined just in terms of function values at the vertices of the uniform grid of level $\ell, \,G^{\ell}\,.$ The only constraint needed to maintain the approximation order of the original QI is to use finite differences approximations of order $(k_1,k_2)$ such that $k_1\ge d_1$ and $k_2\ge d_2$ (see \cite{MS12} for some details on these approximations  and \cite{tesiIurino,ACM15} for their extensions to the tensor--product case). In our tests with $(d_1,d_2)=(3,3)$ we used a finite-difference scheme of order $(k_1,k_2)=(3,3)$ and, in order  to simplify the implementation, we have used the same formulas for the inner points in the whole domains, so taking the necessary gridded function values  from a suitably further enlarged domain. The results obtained show that the performances of $Q^a_{\cal H}$ are essentially the same as the ones of $Q_{\cal H}$ -- see Tables \ref{test13}-\ref{test14d}, which also report the results obtained with the tensor--product version of $Q^a_{\cal H}$, denoted by $Q^{a}$. It is worth noting that, in order to obtain with $Q^a_{\cal H}$ essentially the same error as with $Q^a$, we had to consider meshes with refined areas larger than the ones used with $Q_{\cal H}$. This is in accordance with Theorem \ref{locconv}, since replacing the derivatives with finite-differences schemes enlarges the support of the functionals $\lambda_J^{\ell}$ and, as a consequence, the diameter of the set $C$ defined in \eqref{Cbigdef}.

For brevity, in the tables we use the following notation,
\begin{equation*}
e^a := f-Q^a(f),\qquad e^a_{\cal H} := f-Q^a_{\cal H}(f).
\end{equation*}

\begin{figure}[H] %\begin{center}
\centering
\subfigure[]{\includegraphics[trim=0mm 0mm 0mm 0mm,clip,scale=0.41]{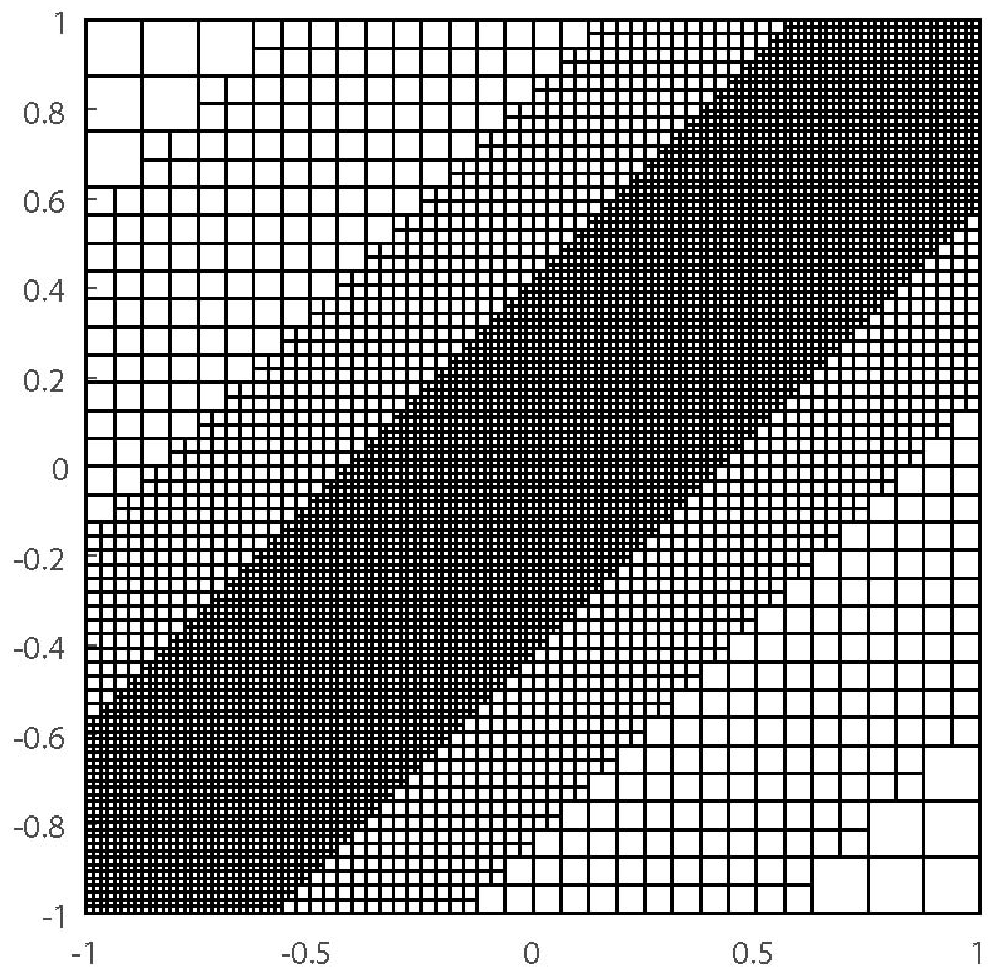}}  %change figure
\subfigure[]{\includegraphics[trim=0mm 0mm 0mm 0mm,clip,scale=0.41]{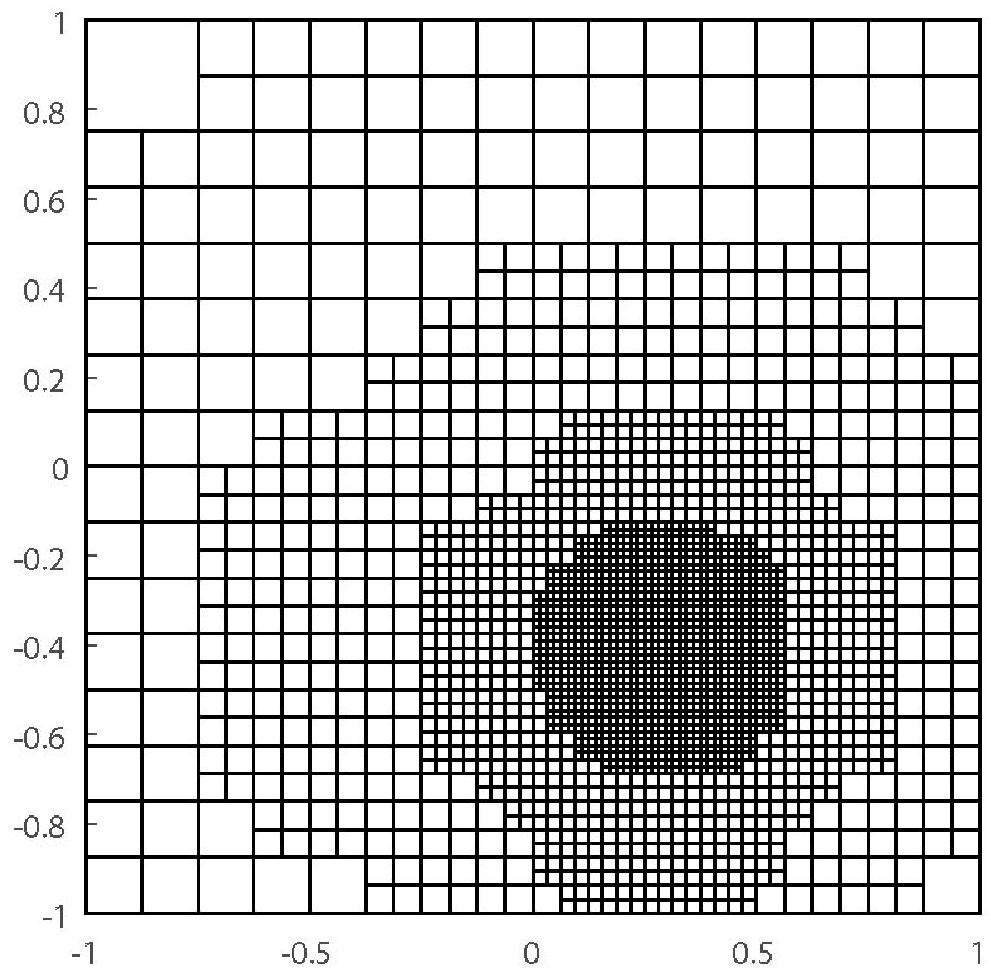}}
\caption{The hierarchical meshes with 5 levels used in the approximation of (a) $f_1$ and (b) $f_2$ with $Q^a_{\cal H}$, $(d_1,d_2)=(3,3)$ and $(k_1,k_2)=(3,3)$.}
%\end{center}
\end{figure}

\begin{table}[H]
\centerline{ 
\begin{footnotesize}\begin{tabular}{SS|SS|SS} \toprule
    {$M$} & {${h}_{M-1} $} & {$\dim(V^{M-1}$)} & {$\Vert e^a\Vert_{ \Omega}$} & {$\dim(S_{\cal H})$} & {$\Vert e^a_{\cal H}\Vert_{ \Omega}$} \\ \midrule
    %1  & {1/2} & {49} & {5.454e-2} & {49} & {5.454e-2} \\
    1  & {1/4} & {121} & {2.538e-2} & {121} & {2.538e-2} \\
    2  & {1/8} & {361} & {4.324e-3} & {361} & {4.324e-3} \\
    3  & {1/16} & {1225} & {4.660e-4} & {1153} & {4.660e-4} \\ 
    4  & {1/32} & {4489} & {2.429e-5} & {3175} & {2.429e-5} \\
    5  & {1/64} & {17161} & {6.323e-7} & {7597} & {6.323e-7} \\ \bottomrule
\end{tabular}\end{footnotesize}
}
\caption{Numerical results for the modified $Q^a_{\cal H}$ (derivatives approximated with order 3 finite-differences scheme) tested on $f_1$ with bi-degree $(3,3)$.}\label{test13}
\end{table}

\begin{table}[H]
\centerline{ 
\begin{footnotesize}\begin{tabular}{S|SSS|SSS} \toprule
    {$M$} & {$\Vert (e_a)_x\Vert_{ \Omega}$} & {$\Vert (e_a)_y\Vert_{ \Omega}$} & {$\Vert (e_a)_{xy}\Vert_{ \Omega}$}& {$\Vert (e^a_{\cal H})_x\Vert_{ \Omega}$} & {$\Vert (e^a_{\cal H})_y\Vert_{ \Omega}$} & {$\Vert (e^a_{\cal H})_{xy}\Vert_{ \Omega}$} \\ \midrule
    %1  & {7.067e-1} & {7.067e-1} &  {6.795e-0} & {7.067e-1} & {7.067e-1} &  {6.795e-0} \\
    1  & {4.302e-1} & {4.302e-1} &  {5.983e-0} & {4.302e-1} & {4.302e-1} &  {5.983e-0} \\
    2  & {1.172e-1} & {1.172e-1} &  {3.057e-0} & {1.172e-1} & {1.172e-1} &  {3.057e-0}\\
    3  & {1.601e-2} & {1.601e-2} &  {6.348e-1} & {1.601e-2} & {1.601e-2} &  {6.348e-1} \\ 
    4  & {1.191e-3} & {1.191e-3} &  {6.229e-2} & {1.191e-3} & {1.191e-3} &  {6.229e-2} \\
    5  & {1.052e-4} & {1.052e-4} &  {6.695e-3} & {1.052e-4} & {1.052e-4} &  {6.695e-3} \\ \bottomrule
\end{tabular}\end{footnotesize}
}
\caption{Errors with respect to the derivatives for $Q^a_{\cal H}$ tested on $f_1$ with $(d_1,d_2)=(3,3)$ and compared with $Q_a$.}\label{test13d}
\end{table}

\begin{table}[H]
\centerline{ 
\begin{footnotesize}\begin{tabular}{SS|SS|SS} \toprule
    {$M$} & {${h}_{M-1} $} & {$\dim(V^{M-1}$)} & {$\Vert e^a\Vert_{ \Omega}$} & {$\dim(S_{\cal H})$} & {$\Vert e^a_{\cal H}\Vert_{ \Omega}$} \\ \midrule
    %1  & {1/2} & {49} & {6.636e-1} & {49} & {6.636e-1} \\
    1  & {1/4} & {121} & {5.043e-1} & {121} & {5.043e-1} \\
    2  & {1/8} & {361} & {1.548e-1} & {352} & {1.548e-1} \\
    3  & {1/16} & {1225} & {2.035e-2} & {772} & {2.035e-2} \\ 
    4  & {1/32} & {4489} & {1.869e-3} & {1417} & {1.869e-3} \\
    5  & {1/64} & {17161} & {8.458e-5} & {2155} & {8.458e-5} \\ \bottomrule
\end{tabular}\end{footnotesize}
}
\caption{Numerical results for the modified $Q^a_{\cal H}$ (derivatives approximated with order 3 finite-differences scheme) tested on $f_2$ with bi-degree $(3,3)$.}\label{test14}
\end{table}

\begin{table}[H]
\centerline{ 
\begin{footnotesize}\begin{tabular}{S|SSS|SSS} \toprule
    {$M$} & {$\Vert (e_a)_x\Vert_{ \Omega}$} & {$\Vert (e_a)_y\Vert_{ \Omega}$} & {$\Vert (e_a)_{xy}\Vert_{ \Omega}$}& {$\Vert (e^a_{\cal H})_x\Vert_{ \Omega}$} & {$\Vert (e^a_{\cal H})_y\Vert_{ \Omega}$} & {$\Vert (e^a_{\cal H})_{xy}\Vert_{ \Omega}$} \\ \midrule
    %1  & {5.719e-0} & {5.712e-0} &  {4.894e+1} & {5.719e-0} & {5.712e-0} &  {4.894e+1} \\
    1  & {5.448e-0} & {5.804e-0} &  {4.937e+1} & {5.448e-0} & {5.804e-0} &  {4.937e+1} \\
    2  & {2.893e-0} & {2.405e-0} &  {2.835e+1} & {2.893e-0} & {2.405e-0} &  {2.835e+1} \\
    3  & {5.013e-1} & {5.085e-1} &  {4.613e-0} & {5.013e-1} & {5.085e-1} &  {4.613e-0} \\ 
    4  & {5.468e-2} & {5.602e-2} &  {5.841e-1} & {5.468e-2} & {5.602e-2} &  {5.841e-1} \\
    5  & {4.195e-3} & {4.124e-3} &  {4.198e-2} & {4.195e-3} & {4.124e-3} &  {9.306e-2} \\ \bottomrule
\end{tabular}\end{footnotesize}
}
\caption{Errors with respect to the derivatives for $Q^a_{\cal H}$ tested on $f_2$ with $(d_1,d_2)=(3,3)$ and compared with $Q_a$.}\label{test14d}
\end{table}

We note that also for the first and second mixed derivatives  the errors produced by  the hierarchical approach are essentially equal to those generated by the corresponding tensor--product QIs.
  
  \section{Conclusions} \label{sec_C}
 The bivariate hierarchical extension  of the BS Hermite spline QI scheme  has been introduced. The convergence properties of the adaptive construction have been presented together with a detailed analysis of performances and costs, also in comparison with a different hierarchical quasi--interpolation method. For the numerical experiments, we constructed and used an automatic mesh refinement algorithm producing hierarchical meshes which in the considered examples allowed to reach the same accuracy of tensor--product quasi--interpolants but considering hierarchical spline spaces with remarkably lower dimension. Furthermore, it has been verified that a variant of the scheme that does not require information on the derivatives has similar features. %The design of cheap automatic strategies for the identification of the suitable subdomain hierarchies  is an additional ingredient to be investigated for the development of a reliable and effective adaptive scheme.

\section*{Acknowledgements}

This work was supported by the programs ``Finanziamento Giovani Ricercatori 2014'' and 
``Progetti di Ricerca 2015'' (Gruppo Nazionale per il Calcolo Scientifico of the Istituto 
Nazionale di Alta Matematica Francesco Severi, GNCS - INdAM) 
and by the project DREAMS (MIUR Futuro in Ricerca RBFR13FBI3).

%%%%%%%%%% % % % % % % % % %  %  %  %  %  %  %  %   %   %   %   %   %    %    %    %    %     %     %     %                                                      

%%%%%%%%%%%%%%%%%%%%%%%%%%%%%%%%%%%%%%%%%%%%%%%%%%%%%%%%%

%%%%%%%%%%%%%%%%%%%%%%%%%%%%%%%%%%%%%%%%%%%%%%%%%%%%%%%%%

\end{document}